\documentclass[conference]{IEEEtran}
\IEEEoverridecommandlockouts
\usepackage{cite}
\usepackage{amsmath,amssymb,amsfonts}

\usepackage{mathtools}
\usepackage{graphicx}
\usepackage{epsf}
\usepackage{mathbbol}
\usepackage{stfloats}
\usepackage{subfigure,color,balance,framed}
\usepackage{verbatim}
\usepackage{textcomp}
\usepackage{xcolor}
\usepackage{comment}
\usepackage{multirow}
\usepackage{makecell}
\usepackage{siunitx}
\usepackage{threeparttable}

\usepackage{booktabs}
\usepackage{balance}
\usepackage{mathbbol}
\usepackage{algorithm}
\usepackage{algorithmicx}
\usepackage{algpseudocode}
\usepackage{mathrsfs}
\usepackage{threeparttable}

\newcommand{\method}[1]{\texttt{#1}}
\newcommand{\tr}{{{\mathsf T}}}
\newcommand{\D}{\textnormal{d}}
\newcommand{\col}{\textnormal{col}}
\newcommand{\des}{\textnormal{des}}
\usepackage{fancyhdr}
\usepackage[colorlinks=true,      % false: boxed links; true: colored links
linkcolor=black,      % color of internal links
citecolor=black,      % color of links to bibliography
filecolor=black,      % color of file links
urlcolor=blue]{hyperref}
\newtheorem{definition}{Definition}

\newtheorem{proposition}{Proposition}
\newtheorem{remark}{Remark}

\newtheorem{problem}{Problem}

\def\BibTeX{{\rm B\kern-.05em{\sc i\kern-.025em b}\kern-.08em
    T\kern-.1667em\lower.7ex\hbox{E}\kern-.125emX}}
\DeclarePairedDelimiter{\norm}{\lVert}{\rVert} 
\begin{document}
\title{Decentralized Robust Data-driven Predictive Control for Smoothing Mixed Traffic Flow}

\author{Xu Shang, Jiawei Wang, and Yang Zheng
\thanks{The work of X. Shang and Y. Zheng is supported by NSF ECCS-2154650 and NSF CMMI-2320697.}% <-this % stops a space
	\thanks{X. Shang and Y. Zheng are with the Department of Electrical and Computer Engineering, University of California San Diego, CA 92093, USA. (x3shang@ucsd.edu; zhengy@ucsd.edu), }%
	\thanks{J. Wang is with the Department of Civil and Environmental Engineering, University of Michigan, Ann Arbor, MI 48109, USA. (jiawe@umich.edu).}%
 }

\maketitle
\thispagestyle{plain}
% \fancyhf[rh]{\thepage}
\begin{abstract}
In a mixed traffic with connected automated~vehicles (CAVs) and human-driven vehicles (HDVs),~data-driven predictive control of CAVs promises system-wide traffic performance improvements. Yet, most existing approaches~focus on a centralized setup, which is computationally unscalable while failing to protect data privacy. The robustness~against~unknown disturbances has not been well addressed either, causing safety concerns. In this paper, we propose a decentralized robust \method{DeeP-LCC} (Data-EnablEd Predictive Leading Cruise Control) approach for CAVs to smooth mixed traffic. In particular, each CAV computes its control input based on locally available data from its involved subsystem. Meanwhile, the interaction between neighboring subsystems is modeled as a bounded disturbance, for which appropriate estimation methods are proposed. Then, we formulate a robust optimization problem and present its tractable computational solutions. Compared with the centralized formulation, our method greatly reduces computation complexity with better safety performance, while naturally preserving data privacy. Extensive traffic simulations validate its wave-dampening ability, safety performance, and computational benefits.
\end{abstract}

\section{Introduction} \label{Introduction}

\IEEEPARstart{U}{ndesirable} traffic waves, such as phantom traffic~jams~\cite{sugiyama2008traffic}, significantly reduce travel efficiency, fuel economy and driving safety. The emergence of connected automated vehicles (CAVs) promises to greatly mitigate this issue~\cite{li2017dynamical}. One typical technology is Cooperative Adaptive Cruise Control (CACC), which groups a series of CAVs~into a platoon and applies cooperative control strategies to mitigate undesired traffic oscillations~\cite{zheng2015stability, milanes2013cooperative}. However, such technologies require a fully CAV environment while the near future will see a transition phase of mixed traffic, where human-driven vehicles (HDVs) coexist with CAVs~\cite{stern2018dissipation, zheng2020smoothing, li2022cooperative}. To enhance their societal benefits, CAVs need to be cooperative with other HDVs that still take the majority in mixed traffic. 

\subsection{CAV Control in Mixed Traffic}

Mixed traffic systems are complex human-in-the-loop cyber-physical systems. Existing research on CAV control in mixed traffic can be divided into model-based and model-free categories, depending on how to address human driving behaviors. Model-based methods typically rely on classical car-following models for HDVs, such as the optimal velocity model (OVM)~\cite{bando1995dynamical} or the intelligent driver model (IDM)~\cite{treiber2000congested}. Lumping the dynamics of HDVs and CAVs, one can derive a parametric model of the entire mixed traffic system~\cite{zheng2020smoothing}, allowing for common model-based synthesis approaches, such as optimal control~\cite{jin2016optimal, wang2021controllability}, $\mathcal{H}_\infty$ control~\cite{di2019cooperative, li2023learning}~and~model predictive control (MPC)~\cite{feng2021robust, guo2021anticipative}. However, it is non-trivial to accurately identify HDVs' car-following behavior, due to the nonlinearity and stochastic nature of human driving behaviors. On the other hand, model-free methods for CAV control, which bypass system identification and directly construct controllers from data, have received increasing attention. For example, reinforcement learning~\cite{wu2021flow, vinitsky2018lagrangian} and adaptive dynamic programming~\cite{gao2016data, huang2020learning} have shown their potential in learning wave-dampening strategies for CAVs. Nonetheless, these model-free methods suffer from several drawbacks, \emph{e.g.}, the lack of interpretability, heavy computation burden~for~offline training, and difficulties of providing safety guarantees.

Recently, a class of data-driven predictive control strategies, which combine learning methods with the well-established MPC, has shown promising results for constrained control~with safety performance~\cite{zhan2022data,lan2021data,wang2023deep}. Notably, a recent study in~\cite{wang2023deep} extends the  Data-EnableEd Predictive Control (\method{DeePC}) technique~\cite{coulson2019data} into Leading Cruise Control~\cite{wang2021leading} (LCC) for mixed traffic, leading to a new notion of \method{DeeP-LCC}. Compared with traditional frameworks such as CACC and Connected Cruise Control~\cite{orosz2016connected}, LCC not only enables the CAVs to adaptively follow the HDVs ahead, but also explores CAVs' potential in actively leading the motion of the HDVs behind. Accordingly, LCC is applicable to general cases of mixed traffic with different CAV formations~\cite{li2022cooperative}, and enables system-wide improvement for global traffic flow~\cite{zheng2020smoothing,wang2021leading}. Note that \method{DeeP-LCC} relies on Willems' Fundamental Lemma~\cite{willems2005note} for data-driven behavior representation of a mixed traffic system. This method requires no prior knowledge of the HDVs' car-following behaviors, and explicitly incorporates input/output constraints (\emph{e.g.}, control saturation and driving safety) into receding horizon control. Both numerical simulations~\cite{wang2023deep} and real-world experiments~\cite{wang2022implementation} have validated the capability of \method{DeeP-LCC} in achieving optimal and safe control for CAVs. 

\subsection{Decentralization and Robustness}

Despite the advantages of data-driven methods in addressing unknown human driving behaviors, the aforementioned studies are still non-trivial in practical traffic systems due to their common centralized setup: a central unit is required to collect trajectory data from all vehicles and to calculate control inputs for all CAVs in a mixed traffic flow (see~\cite{wang2022implementation}~for~implementation details). The centralized setup might be suitable for a small number of vehicles (\emph{e.g.}, one or two CAVs), but is not scalable for real-world traffic due to challenges in real-time computation, data privacy, and communication delays~\cite{gao2016robust,li2017dynamical}. Moreover, the formation captured by the spatial pattern of CAVs, is not static in practical traffic as the surrounding vehicles can join and leave freely~\cite{li2022cooperative}.  In the centralized setup, any changes in the formation require a complete data recollection and controller redesign. To address this,~some recent research has utilized distributed optimization techniques to design distributed data-driven control strategies~\cite{zhan2022data,wang2022distributed}. Particularly, the authors in~\cite{wang2022distributed} solve the original centralized \method{DeeP-LCC} in a distributed manner via the alternating direction method of multiplier (ADMM). However, these distributed algorithms require multiple iterations to compute the control inputs, resulting in high-frequency data exchange between CAVs in a small time period which is often beyond current vehicle-to-vehicle communication capabilities~\cite{lyu2019characterizing}. 

In addition to scalability, robustness against external disturbances also plays a critical role in enhancing the CAVs' potential in practical traffic. Existing research mostly focuses on maximizing the mitigation ratio with respect to external disturbances to improve CAVs' wave-dampening capability. Typical tools include string stability analysis~\cite{zhou2020stabilizing,monteil2018mathcal} or optimal and robust control~\cite{wang2021controllability,mousavi2022synthesis}. % are usually applied. 
However, these studies have not explicitly addressed safety performance against unknown disturbances. Very recently, some safety-related techniques like control barrier function~\cite{zhao2023safety} or tube-based MPC~\cite{feng2021robust} have been proposed, but they require prior knowledge of HDVs' car-following dynamics. In the centralized \method{DeeP-LCC}~\cite{wang2023deep}, a simplistic constant assumption is made for external disturbances, which often fails to accurately reflect human driving behaviors. The mismatch of future external disturbances for \method{DeeP-LCC}~\cite{wang2023deep} may lead to rear-end collisions. While robust data-driven control has been explored in \method{DeePC} for power applications~\cite{huang2021decentralized}, it can not be directly applied to mixed traffic control because of different system dynamics.

\subsection{Contributions}
\begin{figure*}[t]
	\begin{center}
		\includegraphics[width=1\textwidth]{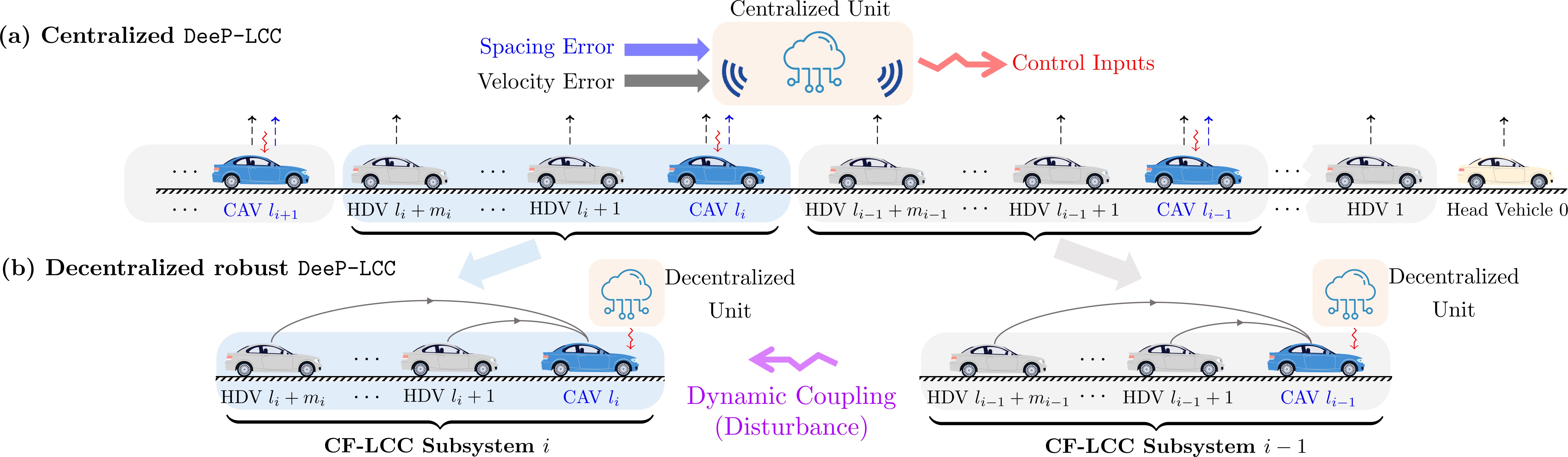}
		\caption{Schematic of centralized and decentralized robsut \method{DeeP-LCC} for CAVs in mixed traffic. (a) Centralized \method{DeeP-LCC}. It collects the data of the whole mixed traffic system, including velocity errors of all vehicles (black dashed arrow) and spacing errors of CAVs collectively (blue dashed arrow), and utilizes these data to design control input of CAVs (red squiggle arrow). (b) Decentralized robust \method{DeeP-LCC}. It decomposes the mixed traffic system into multiple CF-LCC subsystems and only requires locally available data. The dynamic coupling between subsystems is modeled as the disturbance (purple squiggle~arrow).  }
  \label{fig:MixTrafSys}
	\end{center}
 \vspace{-5mm}
\end{figure*}

In this paper, we develop a  decentralized robust data-driven predictive control approach to improve the scalability, robustness, and data privacy of CAV control in mixed traffic (particularly in the spirit of \method{DeeP-LCC}~\cite{wang2023deep}). As shown in Fig. \ref{fig:MixTrafSys}, we first decentralize the problem by decomposing the entire mixed traffic system into multiple CF-LCC (Car-Following Leading Cruise Control) subsystems. Each CAV utilizes the measurable data from its own subsystem to compute safe and optimal control inputs.  We form a decentralized robust optimization problem for online predictive control, in which different disturbance estimation methods are incorporated. We finally provide an efficient computational method to solve the robust optimization problem. Some preliminary results are reported in~\cite{shang2023smoothing}. Our contributions of this paper are as follows: 

\begin{itemize}
    \item We propose a decentralized \method{DeeP-LCC} formulation for CAVs in mixed traffic with locally measured trajectory data of each CF-LCC subsystem. 
 Compared with centralized \method{DeeP-LCC}~\cite{wang2023deep}, it significantly reduces the computation time with a smaller amount of required data. Also, unlike existing distributed optimization approaches~\cite{wang2022distributed,zhan2022data}, our formulation bypasses frequent inter-vehicle communications between neighboring subsystems and naturally contributes to better data privacy. 
 
\item We robustify the decentralized \method{DeeP-LCC} with appropriate disturbance estimation methods to handle the coupling constraints between CF-LCC subsystems. The disturbance bound is estimated based on the system dynamics of the mixed traffic system. Integrating a time-varying disturbance estimation with the aforementioned robustification, our method provides better safety guarantees for each CAV while maintaining comparable wave-damping performance with respect to centralized \method{DeeP-LCC}~\cite{wang2023deep}.

\item Finally, we provide an efficient computational method to solve the decentralized robust \method{DeeP-LCC} problem. The robustification and safety requirements of the mixed traffic system significantly increase the computational complexity of online predictive control. For this aspect, we compare and analyze two computational methods, which are motivated and adapted from robust optimization literature~\cite{bertsimas2011theory, lofberg2012automatic}. We have further developed an efficient automatic transformation that reformulates each decentralized robust \method{DeeP-LCC} problem into its standard conic form, reducing the modeling time compared to general tools, such as YALMIP~\cite{Lofberg2004}. 
\end{itemize}

 Numerical experiments validate the superior performance of our decentralized robust \method{DeeP-LCC} in terms of both computational efficiency and safety performance. In particular, its computation time is $0.058\, \mathrm{s}$, which is $\textbf{85.1\%}$ less than the centralized~\method{DeeP-LCC}, facilitating its real-time computation. Moreover, in our experiments, it achieves $\textbf{0\%}$ violation rate for safety constraints when using both small ($T=700$) and large ($T=1500$) data sets, while these numbers are $\textbf{99\%}$ and $\textbf{88\%}$ for the normal \method{DeeP-LCC} without robustification~\cite{wang2023deep}.

\subsection{Paper Structure}
The remainder of this paper is structured as follows.~Section~\ref{Review of Centralized DeeP-LCC} reviews the basics of centralized \method{DeeP-LCC}.  The formulation of decentralized robust \method{DeeP-LCC} is introduced in Section \ref{Decentralized DeeP-LCC}, and the disturbance estimation methods are discussed in Section \ref{sec:dist_est}. Section \ref{robust_computation} analyzes different computational approaches. Section \ref{Results} illustrates our numerical traffic simulations. Finally, Section \ref{Conclusion} concludes this paper.

\textit{Notations:} We denote $\mathbb{N}$ as the set of natural numbers. We use $\mathbb{0}$ and $I$ as a zero matrix and an identity matrix with a compatible size. For a series of vectors $a_1, a_2, \ldots, a_n$ and matrices $A_1, A_2, \ldots, A_n$ with the same column size, we denote col$(a_1, a_2, \ldots, a_n) = [a_1^\tr, a_2^\tr, \ldots, a_n^\tr]^\tr$ and col$(A_1, A_2, \ldots, A_n) = [A_1^\tr, A_2^\tr, \ldots, A_n^\tr]^\tr$. For a vector $a$ and a square matrix $X$, 
We also denote $\|a\|_X^2 = a^\tr X a$, and diag($D_1, D_2, \ldots, D_n$) denotes a block-diagonal matrix with diagonal blocks $D_1, D_2, \ldots, D_n$. The floor function $\lfloor \alpha \rfloor$ gives the largest integer that is less than $\alpha$. Finally, we denote the Kronecker product between matrices $A$ and $B$ as $A \otimes B$. 

\vspace{2em}
\section{Review of Centralized \method{DeeP-LCC} \\and Problem Statement}
\label{Review of Centralized DeeP-LCC}

In this section, we briefly review the modeling process for the mixed traffic system consisting of multiple CAVs and HDVs and the centralized \method{DeeP-LCC}~\cite{wang2023deep}. In the end, we present the problem statement for this work.

\subsection{Input/Output of Mixed Traffic System}
\label{subsec:IO-mixed}

As shown in Fig.~\ref{fig:MixTrafSys}, we consider a general single-lane mixed traffic system with one head vehicle, indexed as $0$, and $n$ following vehicles, indexed as $1,\ldots,n$ from front to end. Among the following vehicles, there exist $q$ CAVs, $n-q$ HDVs and the $i$-th CAV is indexed as $l_i$ ($1 \leq l_1 < l_2 < \cdots < l_q \leq n$). The number of the HDVs between CAVs $i$ and $i+1$ is represented as $m_i$, \emph{i.e.}, $m_i = l_{i+1}-l_i-1$. %, and we have $l_i+m_i=l_{i+1}-1$. 
We denote the index set $\Omega$ for all the following vehicles, $\mathcal{S}$ for all the CAVs, $\mathcal{F}_i$ for the HDVs between CAV $i$ and CAV $i+1$, and $\mathcal{F}$ for all the HDVs: 
\begin{equation*}
    \begin{aligned}
    \Omega  &= \{1,2,\ldots,n\},  \quad \mathcal{S} =\{l_1,l_2,\ldots,l_q\},\\
    \mathcal{F}_i &= \{l_i+1, \ldots, l_i+m_i\},  \quad \mathcal{F} = \mathcal{F}_1 \cup \ldots \mathcal{F}_n.
    \end{aligned}
\end{equation*}

The position, velocity, and acceleration of vehicle $i$ at time $t$ are denoted as $p_i(t)$, $v_i(t)$ and $a_i(t)$, $i \in \Omega$, respectively. The spacing between vehicle $i$ and its preceding vehicle $i-1$ is denoted as $s_i(t) = p_{i-1}(t) - p_i(t)$, and their relative velocity is $\dot{s}_i (t) = v_{i-1}(t) - v_i(t)$. At an equilibrium state, all vehicles move at an identical velocity $v^*$ and maintain the corresponding equilibrium space $s^*_i$, which may vary from different vehicles.  
We consider the error state of each vehicle with respect to its equilibrium state to facilitate system analysis and controller design. In particular, the velocity error and spacing error for each vehicle are defined as: 
\[
\tilde{v}_i(t) = v_i(t) - v^*, \quad \tilde{s}_i(t) = s_i(t) - s^*_i, \quad i\in\Omega.
\]

We then utilize the aggregation of the velocity and spacing errors of 
the following vehicles to represent the global state $x \in \mathbb{R}^{2n}$ of the mixed traffic system, given by 
\begin{equation} \label{eqn:state}
x(t) = \begin{bmatrix}\tilde{s}_1(t), \tilde{v}_1(t), \tilde{s}_2(t), \tilde{v}_2(t), 
\ldots, \tilde{s}_{n}(t), \tilde{v}_{n}(t)\end{bmatrix}^{\tr}.
\end{equation}
The state \eqref{eqn:state} is not fully measurable. Indeed, it is impractical to obtain the equilibrium spacing $s^*_i, i\in \mathcal{F}$ for HDVs, due to~the~unknown car-following behaviors of HDVs. The measurable output of the traffic system only consists of velocity errors for all vehicles and spacing errors for CAVs, as we can estimate the equilibrium velocity $v^*$ from the head vehicle and the equilibrium spacing $s^*_i$ for CAVs ($i\in\mathcal{S}$) is a designed parameter. Accordingly, the output $y(t) \in \mathbb{R}^{n+q}$ is defined as
\begin{equation} \label{eqn:output}
y(t) = \begin{bmatrix}\tilde{v}_1(t), \tilde{v}_2(t), \ldots, \tilde{v}_{n}(t), \tilde{s}_{l_1}(t), \ldots, \tilde{s}_{l_q}(t)\end{bmatrix}^\tr.
\end{equation}

As widely used before~\cite{zheng2020smoothing,orosz2016connected,wang2021controllability}, we assume the acceleration of each CAV can be directly controlled in this study. We denote $u_i(t)= a_i(t)\in \mathbb{R}$ as the control input of each CAV ($i\in\mathcal{S}$). Then, The control input of the entire mixed traffic system $u(t) \in \mathbb{R}^q$ is defined as 
\begin{equation} \label{eqn:input}
u(t) = \begin{bmatrix}u_{l_1}(t), u_{l_2}(t), \ldots, u_{l_q}(t)\end{bmatrix}^\tr.
\end{equation}
Additionally, the head vehicle's velocity error is considered as an external disturbance signal since it cannot be controlled, which is denoted~as  
\begin{equation} \label{eqn:disturbance}
\epsilon(t) = \tilde{v}_0 = v_0(t) - v^*.
\end{equation}

Upon defining the state~\eqref{eqn:state}, output~\eqref{eqn:output}, input~\eqref{eqn:input} and disturbance~\eqref{eqn:disturbance}, after linearization and discretization, the parametric model of the mixed traffic system can be derived as
\begin{equation}
\label{eqn:ModelMixTraffic}
\left\{
\begin{aligned}
x(k+1) & = A x(k) + B u(k) + H \epsilon(k), \\
y(k) &= C x(k),
\end{aligned}
\right.
\end{equation}
where $k$ denotes the discrete time step and the specific form of the matrices $A \in \mathbb{R}^{2n \times 2n}$, $B\in \mathbb{R}^{2n \times q}, C\in\mathbb{R}^{(n+q)\times2n}$, and $H\in\mathbb{R}^{2n}$ can be found in~\cite[Section II-C]{wang2023deep}; also see~\cite{wang2021leading}.

One goal of CAV control is to design suitable input $u(t)$ in \eqref{eqn:input} to smooth the mixed traffic flow, represented by $x(t)$ in \eqref{eqn:state}, based on the observation $y(t)$ in \eqref{eqn:output}, in the presence of the external disturbance $\epsilon(t)$. We remark that the parametric model \eqref{eqn:ModelMixTraffic}, characterized by matrices  $A, B, C, H$, is not available for controller design due to the unknown HDVs' behavior. To address this, the recently emerging data-driven methods, particularly the centralized \method{DeeP-LCC}~\cite{wang2023deep},  bypass system identification and directly use the input/output trajectories for traffic behavior prediction and CAV control.

\subsection{Basics of Centralized \method{DeeP-LCC}}
The centralized \method{DeeP-LCC}~\cite{wang2023deep} adapts and extends the standard \method{DeePC}~\cite{coulson2019data} for CAV control in mixed traffic. It utilizes a pre-collected input/output trajectory to construct a data-driven behavior representation for the mixed traffic system~\eqref{eqn:ModelMixTraffic}. The pre-collected data needs to be rich enough, which is required by~a~persistent excitation condition.
\begin{definition}[Persistently Exciting]
\label{def:PE}
    The sequence of signal $\omega = \textrm{col}(\omega(1),\omega(2), \ldots, \omega(T))$ with length $T$ ($T \in \mathbb{N}$) is persistently exciting of order $L$ ($L < T$) if its associated Hankel matrix with depth $L$ has full row rank:
    \[\mathcal{H}_L(\omega) = \begin{bmatrix}
        \omega(1) & \omega(2) & \cdots & \omega(T-L+1) \\
        \omega(2) & \omega(3) & \cdots & \omega(T-L+2) \\
        \vdots    & \vdots    & \ddots & \vdots \\
        \omega(L) & \omega(L+1) & \cdots & \omega(T)
    \end{bmatrix}.\]
\end{definition}

We start by collecting a length-$T$ input/output data sequence from the mixed traffic system: 
\begin{equation} \label{eq:offline-data-mixed-traffic}
\begin{aligned}
u^\D & = \textrm{col}(u^\D(1), u^\D(2), \ldots, u^\D(T)) \in \mathbb{R}^{qT}, \\
\epsilon^\D & = \textrm{col}(\epsilon^\D(1), \epsilon^\D(2), \ldots, \epsilon^\D(T)) \in \mathbb{R}^T, \\
y^\D & = \textrm{col}(y^\D(1), y^\D(2),\ldots,y^\D(T)) \in \mathbb{R}^{(n+q)T},
\end{aligned}
\end{equation}
and then partition associated Hankel matrices into past data with length $T_\textrm{ini}$ and future data with length $N$ 
\begin{equation}
\label{eqn:Hankel-partitioned}
\begin{bmatrix}
    U_{\textrm{P}} \\
    U_{\textrm{F}} 
\end{bmatrix} := \mathcal{H}_L(u^\textrm{d}), \ 
\begin{bmatrix}
    E_{\textrm{P}} \\
    E_{\textrm{F}} 
\end{bmatrix} := \mathcal{H}_L(\epsilon^\textrm{d}), \
\begin{bmatrix}
    Y_{\textrm{P}} \\
    Y_{\textrm{F}} 
\end{bmatrix} := \mathcal{H}_L(y^\textrm{d}),
\end{equation}
where $L = T_\textrm{ini}+N$, and $U_{\textnormal{P}}$ and $U_{\textnormal{F}}$ consist of the first $T_\textnormal{ini}$ rows and the last $N$ rows of $\mathcal{H}_L(u^\textnormal{d})$, respectively (similarly for $E_{\textrm{P}}$ and $E_{\textrm{F}}$, $Y_{\textrm{P}}$ and $Y_{\textrm{F}}$). We note that the mixed traffic system \eqref{eqn:ModelMixTraffic} is controllable~\cite{wang2021leading} and we have the following result.
\begin{proposition}[\!\!{\cite[Proposition 2]{wang2023deep}}] \label{proposition:data-rep-central}
Let the most recent past input trajectory with length $T_\textrm{ini}$ and the future input trajectory with length $N$ as $u_\textrm{ini} = \textrm{col}(u(t-T_\textrm{ini}),u(t-T_\textrm{ini}+1),\ldots, u(t-1))$ and $u = \textrm{col}(u(t), u(t+1),\ldots, u(t+N-1))$, respectively (similarly for $\epsilon_\textrm{ini}, \epsilon, y_\textrm{ini}, y$). If the offline data $\bar{u}^\D := \textrm{col}(u^\D, \epsilon^\D)$ is persistently exciting of order $L+2n$, then the online data $\textrm{col}(u_\textrm{ini}, \epsilon_\textrm{ini}, y_\textrm{ini}, u, \epsilon, y)$ is a valid length-$L$ trajectory of \eqref{eqn:ModelMixTraffic} if and only if there exists $g \in \mathbb{R}^{T-T_\textrm{ini} -N+1}$ such that
\begin{equation}
\label{eqn:predictor}
\begin{bmatrix}
    U_\textrm{P} \\
    E_\textrm{P} \\
    Y_\textrm{P} \\
    U_\textrm{F} \\
    E_\textrm{F} \\
    Y_\textrm{F}
\end{bmatrix} g
=
\begin{bmatrix}
    u_\textrm{ini} \\
    \epsilon_\textrm{ini} \\
    y_\textrm{ini} \\
    u \\
    \epsilon \\
    y
\end{bmatrix}.
\end{equation}
If $T_{\mathrm{ini}} \ge l$, where $l$ denotes the lag of system \eqref{eqn:ModelMixTraffic}, $y$ is uniquely determined from~\eqref{eqn:predictor} for each pair of $u_\textrm{ini}$, $\epsilon_\textrm{ini}$, $y_\textrm{ini}$, $u$ and $\epsilon$.
    
\end{proposition}

The data-driven behavior representation~\eqref{eqn:predictor} holds for the linear time-invariant mixed traffic system~\eqref{eqn:ModelMixTraffic}. Given the nonlinear and nondeterministic nature of practical mixed traffic, the centralized \method{DeeP-LCC}~\cite{wang2023deep} further adapts~\eqref{eqn:predictor} and solves  the following optimization problem at each time step $t$:
\begin{subequations}
\label{eqn:cDeeP-LCC}
\begin{align}
\min_{g, \sigma_y, u, y}  \quad & J(y, u, g, \sigma_y) \label{eqn:cDeeP-LCC-a}
\\
\textrm{subject~to} \quad & 
\begin{bmatrix}
U_\textnormal{P}\\
E_\textnormal{P}\\
Y_\textnormal{P}\\
U_\textnormal{F}\\
E_\textnormal{F}\\
Y_\textnormal{F}\\
\end{bmatrix} g
= 
\begin{bmatrix}
u_\textnormal{ini}\\
\epsilon_\textnormal{ini}\\
y_\textnormal{ini}\\
u\\
\epsilon\\
y\\
\end{bmatrix} + 
\begin{bmatrix}
0\\
0\\
\sigma_y\\
0\\
0\\
0\\
\end{bmatrix}, \label{eqn:cDeeP-LCC-b}
\\
& \epsilon = \mathbb{0}_N, \quad y \in \mathcal{Y}, \quad u \in \mathcal{U}, \label{eqn:cDeeP-LCC-c}
\end{align}
\end{subequations} 
where $(u_\textnormal{ini}, \epsilon_\textnormal{ini}, y_\textnormal{ini})$ is an initial trajectory that needs to be updated at each time step. We will detail the cost function $J(y,u,g,\sigma_y)$, slack variable $\sigma_y$, output constraint $\mathcal{Y}$ and input constraint $\mathcal{U}$ in Section \ref{subsec:rob_formu} (see \eqref{eqn:objFun}, \eqref{eqn:safety}, \eqref{eqn:inputlimit}). One key aspect of \method{DeeP-LCC} \eqref{eqn:cDeeP-LCC} is the data-driven behavior predictor \eqref{eqn:cDeeP-LCC-b} from Proposition \ref{proposition:data-rep-central}, which requires no model information. For the future disturbance $\epsilon$, its value is estimated as $\mathbb{0}_N$ in \cite{wang2023deep} (see \eqref{eqn:cDeeP-LCC-c}) based on the assumption that the head vehicle always tries to maintain its equilibrium state.

\subsection{Problem Statement} \label{subsection:problem-statement}

Despite its promising numerical~\cite{wang2023deep} and experimental~performance~\cite{wang2022implementation}, the original \method{DeeP-LCC}~\eqref{eqn:cDeeP-LCC} has two major~limits:  1) its centralized control setup and 2) the simplistic zero assumption of the future external disturbance in \eqref{eqn:cDeeP-LCC-c}.  %, \emph{i.e.}, the velocity error for the head vehicle. 
In particular, as shown in Fig. \ref{fig:MixTrafSys}(a), the centralized \method{DeeP-LCC} requires global information of the whole traffic system to solve~\eqref{eqn:cDeeP-LCC}. Given a practical mixed traffic system with multiple CAVs and HDVs, the length of pre-collected data~\eqref{eq:offline-data-mixed-traffic} might be huge, affecting data privacy, communication burden, and numerical efficiency of solving~\eqref{eqn:cDeeP-LCC} in real time. Moreover, the zero future velocity error estimation (\emph{i.e.}, $\epsilon = \mathbb{0}_N$) is not consistent with the real-world driving behaviors, and strong oscillations may happen during the occurrence of traffic waves. An inaccurate estimation results in a mismatch between prediction and real traffic behavior, posing safety concerns (rear-end collisions). 

In this work, we aim to design controllers for each CAV individually based on its local data to track the velocity of the head vehicle (or a given reference velocity) and ensure the CAV keeps a desired equilibrium spacing. Our first objective is to address the centralized and zero disturbance settings. %, formally stated as follows:
\begin{problem}[Decentralization and Robustification]
    Develop a decentralized and robust formulation of \method{DeeP-LCC}, which relies on locally available data and incorporates explicit considerations of potential future disturbances. %
\end{problem}

Our strategy is to decompose the entire mixed traffic system into a series of subsystems, and formulate a robust optimization problem against an uncertainty set of future disturbances, which will be detailed in Section~\ref{Decentralized DeeP-LCC}. Then, we need an appropriate uncertainty estimation for the future disturbance~$\epsilon$. This is our second objective, addressed in Section \ref{sec:dist_est}. 
\begin{problem}[Disturbance Estimation and Approximation]
    Estimate an accurate disturbance set from the historical data of the preceding vehicle's velocity error and develop a low-dimension approximation.
\end{problem}

Finally, the third objective is to establish an efficient computation method for the robust data-driven predictive control. %, which will be discussed in Section \ref{robust_computation}.
\begin{problem}[Efficient Computations]
    Develop a computationally efficient method to solve the decentralized robust \method{DeeP-LCC} online at each time step.
\end{problem}

For this objective, we will exploit the problem structure of \method{DeeP-LCC} and adapt techniques from robust optimization~\cite{bertsimas2011theory}~\cite{lofberg2012automatic}. Details will be presented in Section \ref{robust_computation}.

\section{Decentralized Robust \method{DeeP-LCC}} \label{Decentralized DeeP-LCC}
In this section, we present a new decentralized robust \method{DeeP-LCC} formulation to control CAVs in mixed traffic. In particular, we first divide the whole traffic system in Fig.~\ref{fig:MixTrafSys} into multiple subsystems, called  Car-Following LCC (CF-LCC). We then derive a data-driven representation for each CF-LCC subsystem, allowing for its local \method{DeeP-LCC} formulation.  

\subsection{Input/Output of CF-LCC Subsystem}
As shown in Fig. \ref{fig:MixTrafSys}, by exploiting the cascading structure, the entire mixed traffic with $q$ CAVs can be naturally divided into $q$ subsystems, which are called CF-LCC in~\cite{wang2021leading}. Each of the CF-LCC subsystems represents a small mixed traffic system with one CAV $l_i$ and its following $m_i$ HDVs (denoted by $\mathcal{F}_i$). The CAV leads the motion of the HDVs behind it while following one single vehicle immediately ahead, which is typically HDV $l_i-1$ (but~could also be CAV $l_{i-1}$ if $\mathcal{F}_{i-1} =~\varnothing$). Although the centralized \method{DeeP-LCC} may be directly applied to CF-LCC subsystems, a naive application can lead to deteriorated control performance and pose safety concerns (see Section \ref{Results}). 

Here, we focus on the dynamics of each CF-LCC subsystem and \textit{view the coupling dynamics between two neighboring CF-LCC subsystems as an external disturbance}. Similar to~\eqref{eqn:state}-\eqref{eqn:disturbance} for the global mixed traffic system, the state $x_i(t) \in \mathbb{R}^{2(m_i+1)}$ of the $i$-th CF-LCC is defined as
\begin{equation*}
    x_i(t) \!= \! [\tilde{s}_{l_i}\!(t),\! \tilde{v}_{l_i}\!(t),\! \tilde{s}_{l_i+1}\!(t),\! \tilde{v}_{l_i+1}\!(t),\!
     \cdots\!,\! \tilde{s}_{l_i+m_i}\!(t),\! \tilde{v}_{l_i+m_i}\!(t)]^\tr,
\end{equation*}
and the output $y_i(t) \in \mathbb{R}^{m_i+2}$ is defined as
\[
y_i(t) = [\tilde{v}_{l_i}(t), \tilde{v}_{l_i+1}(t), \tilde{v}_{l_i+2}(t), \ldots, \tilde{v}_{l_i+m_i}(t), \tilde{s}_{l_i}(t)]^\tr,
\]
and the control input $u_i(t)$ is the control signal of the  CAV, which is $u_{l_i}(t) \in \mathbb{R}$. In addition, the velocity error of the preceding vehicle is considered as an external disturbance:  
\begin{equation} \label{eq:velocity-error}
\epsilon_i(t) = \tilde{v}_{l_i-1} = v_{l_i-1}-v^* \in \mathbb{R}.
\end{equation}

Similar to~\eqref{eqn:ModelMixTraffic}, we can write the linearized discrete-time state-space model of the $i$-th CF-LCC subsystem as 
\begin{equation}
\label{eqn:ModelCFLCC}
\left\{
\begin{aligned}
x_i(k+1) & = A_i x_i(k) + B_i u_i(k) + H_i \epsilon_i(k), \\
y_i(k) &= C_i x_i(k),
\end{aligned}
\right.
\end{equation}  
where $A_i, B_i, C_i$ and $H_i$ can be found in~\cite{wang2021leading}. The controllability and observability of the CF-LCC subsystem under some mild conditions are also proved in~\cite{wang2021leading}, which enable the design of the optimal controller and the estimation of the state from the output.

\subsection{Data-driven Representation and Local \method{DeeP-LCC}}
\label{subsec:rob_formu}

The CF-LCC model~\eqref{eqn:ModelCFLCC} (\emph{i.e.}, the knowledge of $A_i, B_i, C_i$ and $H_i$) is unknown in practice. We here derive a data-driven representation for the dynamics of each CF-LCC subsystem, allowing for the construction of the local \method{DeeP-LCC}. 

\vspace{0.2em}
\textbf{Offline Data Collection:}  We collect an input/output trajectory of length $T_i$ for the $i$-th subsystem: 
\[
\begin{aligned}
u_i^\textnormal{d} &= \textrm{col}(u_i^\textnormal{d}(1), u_i^\textnormal{d}(2),\ldots,u_i^\textnormal{d}(T_i))\in \mathbb{R}^{T_i}, \\
\epsilon_i^\textnormal{d} &= \textrm{col}(\epsilon_i^\textnormal{d}(1), \epsilon_i^\textnormal{d}(2),\ldots,\epsilon_i^\textnormal{d}(T_i))\in \mathbb{R}^{T_i},\\
y_i^\textnormal{d} &= \textrm{col}(y_i^\textnormal{d}(1), y_i^\textnormal{d}(2),\ldots,y_i^\textnormal{d}(T_i))\in \mathbb{R}^{(m_i+2)T_i}.
\end{aligned}
\] 
These collected data are then used to construct the Hankel matrix of order $L$, which is partitioned as follows: 
\begin{equation} \label{eq:Hankel-CF-LCC}
\begin{bmatrix}
    U_{i,\textnormal{P}} \\
    U_{i,\textnormal{F}} 
\end{bmatrix} \!\! := \!\mathcal{H}_L(u_i^\textnormal{d}),\; 
\begin{bmatrix}
    E_{i,\textnormal{P}} \\
    E_{i,\textnormal{F}} 
\end{bmatrix} \!\!:= \!\mathcal{H}_L(\epsilon_i^\textnormal{d}),\; 
\begin{bmatrix}
    Y_{i,\textnormal{P}} \\
    Y_{i,\textnormal{F}} 
\end{bmatrix} \!\!:= \!\mathcal{H}_L(y_i^\textnormal{d}),
\end{equation}
where $U_{i,\textnormal{P}}$ and $U_{i,\textnormal{F}}$ consist of the first $T_\textnormal{ini}$ rows and the last $N$ rows of $\mathcal{H}_L(u_i^\textnormal{d})$, respectively (similarly for $E_{i,\textnormal{P}}$ and $E_{i,\textnormal{F}}$, $Y_{i,\textnormal{P}}$ and $Y_{i,\textnormal{F}}$).

\vspace{0.2em}
\textbf{Online Behavior Representation:} 
We will utilize the Hankel matrix \eqref{eq:Hankel-CF-LCC} for online behavior predictions of CF-LCC. 
We have the following result.

\begin{proposition}[\!\!{\cite[Lemma 2]{wang2022distributed}}] \label{proposition:data-rep-decen}
At time step $k$, denote the most recent past input sequence $u_{i, \textnormal{ini}}$ with length $T_{\textnormal{ini}}$ and the future input sequence $u_i$ with length $N$ as 
\begin{align*}
 u_{i, \textnormal{ini}} &= \textrm{col}(u_i(k-T_\textnormal{ini}), u_i(k-T_\textnormal{ini}+1),\ldots,u_i(k-1)), \\
 u_i &= \textrm{col}(u_i(k), u_i(k+1), \ldots , u_i(k+N-1)),
\end{align*}
and $\epsilon_{i, \textnormal{ini}}$, $\epsilon_i$, $y_{i,\textnormal{ini}}$ and $y_i$ are denoted similarly.  %follows similar definitions. 
If the data sequence $\bar{u}_i^\D = \textrm{col}(u_i^\textnormal{d}, \epsilon_i^\textnormal{d})$ is persistently exciting of order $L+2(m_i+1)$ (where $L = T_\textnormal{ini} + N$), then 
the sequence $\textrm{col}(u_{i,\textnormal{ini}}, \epsilon_{i,\textnormal{ini}},$ $ y_{i,\textnormal{ini}},u_i,\epsilon_i,y_i)$ is a valid trajectory with length $L$ of the CF-LCC subsystem~\eqref{eqn:ModelCFLCC},  if and only if there exists a vector $g_i \in \mathbb{R}^{T_i-L+1}$ such that 
\begin{equation} \label{eq:CF-LCC-local-representation}
\begin{bmatrix}
U_{i,\textnormal{P}}\\
E_{i,\textnormal{P}}\\
Y_{i,\textnormal{P}}\\
U_{i,\textnormal{F}}\\
E_{i,\textnormal{F}}\\
Y_{i,\textnormal{F}}\\
\end{bmatrix} g_i
= 
\begin{bmatrix}
u_{i,\textnormal{ini}}\\
\epsilon_{i,\textnormal{ini}}\\
y_{i,\textnormal{ini}}\\
u_i\\
\epsilon_i\\
y_i\\
\end{bmatrix}.
\end{equation}
If $T_\textnormal{ini}$ is no smaller than the lag of system \eqref{eqn:ModelCFLCC}, then the future output $y_i$ is unique for fixed $(u_{i,\textnormal{ini}}, \epsilon_{i,\textnormal{ini}}, y_{i,\textnormal{ini}}, u_i, \epsilon_i)$.
\end{proposition}

This proposition establishes a data-driven representation of the dynamics of the CF-LCC subsystem \eqref{eqn:ModelCFLCC}, \emph{i.e.}, all of its valid trajectories can be expressed as a linear combination of the pre-collected trajectories, which are encoded in the Hankel matrices~\eqref{eq:Hankel-CF-LCC}. Accordingly, one can directly utilize the local pre-collected trajectories $(u_i^\textnormal{d}, \epsilon_i^\textnormal{d}, y_i^\textnormal{d})$ to predict the future output $y_i$ given the input $u_i$, the disturbance $\epsilon_i$, and  the initial condition $(u_{i,\textnormal{ini}}, \epsilon_{i,\textnormal{ini}}, y_{i,\textnormal{ini}})$.

Using local data representation \eqref{eq:CF-LCC-local-representation}, we formulate the local $\method{DeeP-LCC}$ for the $i$-th CF-LCC at each time step $k$ as 
\begin{subequations}\label{eqn:localDeeP}
\begin{align}
\min_{g_i, \sigma_{y_i}, u_i, \epsilon_i, y_i}  \quad & V_i(u_i, y_i) + \lambda_{g_i} ||g_i||_2^2 + \lambda_{y_i} ||\sigma_{y_i}||_2^2 \label{eqn:objFun} \\
\textrm{subject~to} \quad & 
\begin{bmatrix}
U_{i,\textnormal{P}}\\
E_{i,\textnormal{P}}\\
Y_{i,\textnormal{P}}\\
U_{i,\textnormal{F}}\\
E_{i,\textnormal{F}}\\
Y_{i,\textnormal{F}}\\
\end{bmatrix} g_i
= 
\begin{bmatrix}
u_{i,\textnormal{ini}}\\
\epsilon_{i,\textnormal{ini}}\\
y_{i,\textnormal{ini}}\\
u_i\\
\epsilon_i\\
y_i\\
\end{bmatrix} + 
\begin{bmatrix}
0\\
0\\
\sigma_{y_i}\\
0\\
0\\
0\\
\end{bmatrix}\label{eqn:equality},\\
& \tilde{s}_{\min} \le G_1 y_i\le \Tilde{s}_{\max}, \label{eqn:safety}\\
& u_{\min} \le u_i \le u_{\max}, \label{eqn:inputlimit}\\
& \epsilon_i = \epsilon_{i,\textrm{est}} \label{eqn:estimator},
\end{align}
\end{subequations}
where 
$V_i(u_i, y_i)$ penalizes the output deviation from equilibrium states and the energy of the input: %, defined as 
\begin{equation}
\label{eq:CL_cost}
V_i(u_i, y_i) = ||u_i||_{R_i}^2 + ||y_i||_{Q_i}^2,
\end{equation}
with $R \in \mathbb{S}_{+}^{N\times N}$ and $Q \in \mathbb{S}_+^{N(m_i+2)\times N(m_i+2)}$. % are positive definite matrices. 
The constraint~\eqref{eqn:equality} is used for online behavior prediction (see Proposition~\ref{proposition:data-rep-decen}), while the constraints~\eqref{eqn:safety} and~\eqref{eqn:inputlimit} are used to ensure safety and driving comfort, respectively. Precisely, $G_1 = I_N \otimes \begin{bmatrix}\mathbb{0}_{1\times (m_i+1)},\; 1\end{bmatrix}$ selects the spacing of the CAV from the output vector $y_i$, and $\tilde{s}_{\min},\tilde{s}_{\max}$ and $u_{\max},u_{\min}$ denote the upper and lower bounds for the spacing error and the control input, respectively.  
The constraint~\eqref{eqn:estimator} denotes the estimation of future velocity errors of the preceding vehicle.
 
Note that Proposition~\ref{proposition:data-rep-decen} only works for LTI mixed traffic systems with no measurement noises and the requirement of the persistent excitation is discussed in detail in the Appendix \ref{appendix:PE}. Yet, any practical mixed traffic system is nonlinear with noises, and therefore, we introduce the regularization term $||g_i||_2^2$ to mitigate overfitting caused by noisy data, the extra decision variable $\sigma_y$ to ensure the feasibility of the optimization problem, and the term $\|\sigma_y\|_2^2$ to penalize the deviation from the initial condition (which is consistent with the original \method{DeePC}~\cite{coulson2019data}). Coefficients $\lambda_g$ and $\lambda_y$ are weighting parameters,  and $\lambda_y$ is usually large enough to satisfy the system's initial condition; see~\cite{dorfler2022bridging,shang2023convex} for details. The closed-loop asymptotical stability of the data-driven controller that has a similar form as \eqref{eqn:localDeeP} is proved in~\cite{berberich2020data} for LTI systems, and the recursive feasibility and practical exponential stability is demonstrated in~\cite{berberich2022linear} for nonlinear systems, which provide potential theoretical guarantees for \method{dDeeP-LCC}.

\begin{algorithm}[t] 
	\caption{Decentralized Robust \method{DeeP-LCC}}
	\label{Alg:dDeeP-LCC}
	\begin{algorithmic}[1]
		\Require
		Pre-collected local data $(u_i^{\textrm{d}},\epsilon_i^{\textrm{d}},y_i^{\textrm{d}})$ for $i$-th subsystem, initial time step $k_0$, terminal time step $k_f$;
		\State Construct data Hankel matrices for input, disturbance and output as $U_{i,\textrm{P}} , U_{i,\textrm{F}}, E_{i,\textrm{P}} , E_{i,\textrm{F}}, Y_{i,\textrm{P}}, Y_{i,\textrm{F}}$;
		\State Initialize the most recent past traffic data $(u_{i,\textrm{ini}},\epsilon_{i,\textrm{ini}},y_{i,\textrm{ini}})$ before the initial time $k_0$;
		\While{$k_0 \leq k \leq k_f$}
  \State Estimate $\mathcal{W}_i$ from $\epsilon_{i, \textrm{ini}}$;
		\State Solve~\eqref{eqn:dDeeP} and get optimal predicted input $u_i^*= \textcolor{white}{\quad \ \, }\textrm{col}(u_i^*(k),u_i^*(k+1),\ldots,u_i^*(k+N-1))$;
		\State Apply the input $u_i(k) \leftarrow u_i^*(k)$ to the $i$-th CAV;
		\State $k \leftarrow k+1$;
            \State Update past local data $(u_{i, \textrm{ini}},\epsilon_{i,\textrm{ini}},y_{i,\textrm{ini}})$;
		\EndWhile
	\end{algorithmic}
\end{algorithm}

\begin{remark}[Coupling dynamics and estimation of future velocity errors] \label{remark:coupling-dynamics}
    The local $\method{DeeP-LCC}$ in~\eqref{eqn:localDeeP} designs the CAV's control input based on an estimated future velocity error trajectory $\epsilon_i$ of the preceding vehicle. A naive choice is the simplistic zero assumption in centralized \method{DeeP-LCC} (\emph{i.e.},~$\epsilon_{i,\textrm{est}} = \mathbb{0}_N$). This inaccurate estimation will cause a mismatch between the prediction and real behavior, which degrades control performance and may lead to a collision. In~\cite{wang2022distributed}, the value of $\epsilon_{i,\textrm{est}}$ is approximated through inter-vehicle data exchange via a distributed algorithm, where impractical high-frequency data exchange is required. In this paper, we estimate a set of potential future velocity error trajectories~$\mathcal{W}_i$, which is more likely to contain the real trajectory. We then require the online optimization~\eqref{eqn:localDeeP} to consider the worst performance against all trajectories in $\mathcal{W}_i$ and establish a robust formulation. \hfill $\square$
\end{remark}

\subsection{Decentralized Robust \method{DeeP-LCC} Formulation}

As discussed in Remark \ref{remark:coupling-dynamics}, instead of considering a single velocity error trajectory, we introduce a velocity error set $\mathcal{W}_i$ as the estimation, \emph{i.e.}, $\epsilon_i \in \mathcal{W}_i$ (the design of this set $\mathcal{W}_i$ will be detailed in Section~\ref{sec:dist_est}). Our key idea here is to formulate an online robust optimization problem by considering the worst performance against $\mathcal{W}_i$. The decentralized robust \method{DeeP-LCC} formulation is formally presented as follows:
\begin{equation}\label{eqn:dDeeP}
\begin{aligned}
\min_{g_i, \sigma_{y_i}, u_i, y_i} \ \max_{\epsilon_i \in \mathcal{W}_i} \quad & V_i(u_i, y_i) + \lambda_{g_i} ||g_i||_2^2 + \lambda_{y_i} ||\sigma_{y_i}||_2^2  \\
\textrm{subject~to} \quad & \eqref{eqn:equality}, \ \eqref{eqn:safety}, \ \eqref{eqn:inputlimit}. %\\
\end{aligned}
\end{equation}

For notational simplicity, we will use $\method{cDeeP-LCC}$ and $\method{dDeeP-LCC}$ to represent centralized \method{DeeP-LCC}~\eqref{eqn:cDeeP-LCC} and decentralized robust $\method{DeeP-LCC}$ \eqref{eqn:dDeeP}, respectively. 
Compared with the basic version~\eqref{eqn:localDeeP}, this robust formulation~\eqref{eqn:dDeeP} improves the online behavioral prediction and thus provides a better safety performance (our numerical experiments also validate this). In the implementation of $\method{dDeeP-LCC}$, $\eqref{eqn:dDeeP}$ is solved in a receding horizon manner, and we re-estimate $\mathcal{W}_i$ iteratively based on the updated velocity errors of the preceding vehicle at each time $k$ (see Section \ref{sec:dist_est} for details). Overall, Algorithm \ref{Alg:dDeeP-LCC} lists the procedure of $\method{dDeeP-LCC}$.

\begin{figure}[t]
\vspace{-1mm}
\centering
{\includegraphics[width=0.48\textwidth]{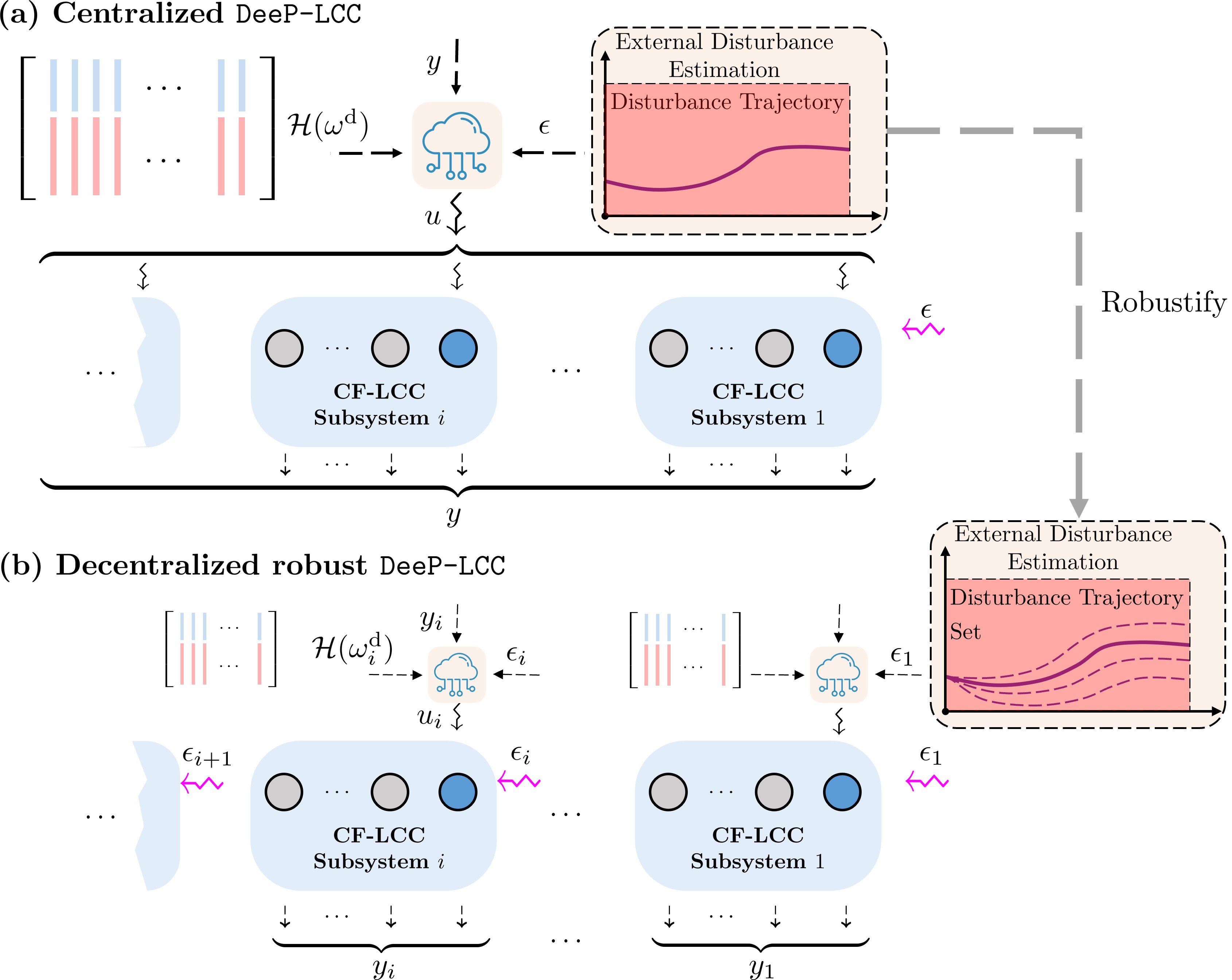}}
\vspace{-1mm}
\caption{Schematic of comparison between centralized and decentralized robust \method{DeeP-LCC}. The Hankel matrix is partitioned into the past trajectories (represented by blue columns) and future trajectories (represented by red columns). The future external disturbance is represented by purple squiggle~arrows.
}
\label{fig:Controller_Comp}
\end{figure}

\begin{remark}[Centralized vs decentralized formulations]
\label{subsec:cen-decen-comp}
In general, \method{cDeeP-LCC} considers the entire mixed traffic system with aggregated input $u(t)$ and output $y(t)$ (see Fig. \ref{fig:Controller_Comp}(a)). In contrast, \method{dDeeP-LCC} focuses on each CF-LCC subsystem (see Fig. \ref{fig:Controller_Comp}(b)) which is more scalable in terms of computation, communication and data privacy. Particularly, the smaller size of the subsystem leads to a lower complexity of the optimization problem, which in turn increases computational efficiency. Also, \method{dDeeP-LCC} requires no communications between subsystems so that the data can be kept inside the local system, which naturally improves data privacy. 
%As each decentralized unit works individually, the failure of one of them shall not affect the functionality of others, while the control of the entire mixed traffic system is messed up if the centralized unit fails. 
Indeed, the CF-LCC subsystems are coupled in a cascading structure (\emph{i.e.}, the velocity error of the preceding vehicle shown in \eqref{eq:velocity-error}, comes from the dynamics of the $(i-1)$-th CF-LCC subsystem), forming the mixed traffic system. Note that the external disturbance $\epsilon$ is simply estimated as zeros in \method{cDeeP-LCC} which is unrealistic and can lead to an emergency if it is directly applied to subsystems (see Section \ref{subsec:safety} for details), while a set of trajectories ($\epsilon_i \in \mathcal{W}_i$) is under consideration in \method{dDeeP-LCC} which improves the safety performance. \hfill $\square$
\end{remark}

\begin{remark}[Data requirement in \method{dDeeP-LCC}] 
The centralized $\method{cDeeP-LCC}$ collects inputs and outputs of all vehicles to construct the associated Hankel matrix for the entire mixed traffic system, while $\method{dDeeP-LCC}$ forms the Hankel matrix for each CF-LCC subsystem with only the local input and output data. Thus, for \method{dDeeP-LCC}, the data can be collected individually for each subsystem and it requires less pre-collected data, which is beneficial for both offline data collection and online behavioral prediction. When the formation pattern of mixed traffic changes due to the free joining and leaving of surrounding vehicles (\emph{i.e.}, the parameters for \eqref{eqn:ModelMixTraffic} can be time-varying), only the data from the affected subsystem needs to be recollected. Furthermore, the smaller system's dimension leads to less complexity for online optimization. In particular, persistently exciting conditions in Propositions \ref{proposition:data-rep-central} and \ref{proposition:data-rep-decen} require Hankel matrices $\mathcal{H}_{L+2n}(\bar{u}^\D)$ and $\mathcal{H}_{L+2(m_i+1)}(\bar{u}_i^\D)$ to have full row rank. To make them square matrices, the minimum data lengths $L_\textrm{cen}$ and $L_\textrm{decen}$ for \method{cDeeP-LCC} and \method{dDeeP-LCC} are, respectively,  as
\begin{subequations}
\label{eqn:data-length}
\begin{align}
    L_\textrm{cen} &= (q+2)(L+2n)-1, \\  L_\textrm{decen} &=3(L+2(m_i+1))-1.
\end{align}
\end{subequations}
Recall that $q$ is the number of CAVs. Since the size of each subsystem is much smaller than the entire system (\emph{i.e.,} we have $q+2 \ge 3$, $n \gg m_i+1$ generally), the $\method{dDeeP-LCC}$ needs much less data for a smaller Hankel matrix; see Fig.~\ref{fig:Controller_Comp}. Thus, $\method{dDeeP-LCC}$ has a smaller data representation in \eqref{eq:CF-LCC-local-representation} which will improve the online computation efficiency. \hfill $\square$
\end{remark}

\section{Disturbance Estimation and Approximation}
\label{sec:dist_est}

Estimating the set of velocity error trajectories ($\epsilon_i \in \mathcal{W}_i$) is a key step in  $\method{dDeeP-LCC}$ that decomposes the whole~system since CF-LCC subsystems are dynamically coupled through velocity perturbation from their preceding vehicles. 
In this section, we first introduce three disturbance estimation approaches based on different assumptions of the preceding vehicle.~We then provide a low-dimensional approximation method via a down-sampling strategy for more efficient implementation.

\subsection{Uncertainty Quantification}
\label{subsec:uncert-quant}
Our general strategy is to bound the set of velocity errors using an $N$-dimensional polytope
\begin{equation} \label{eq:uncertainty-set-W}
\mathcal{W}_i = \{\epsilon_i \in \mathbb{R}^N \ | \ A_\epsilon \epsilon_i \le b_\epsilon \},
\end{equation}
where $A_\epsilon = \textrm{col}(I_N, -I_N)$, $b_\epsilon = \textrm{col}(\epsilon_{i,\max}, \epsilon_{i,\min})$ and $\epsilon_{i,\min}, \epsilon_{i,\max} \in \mathbb{R}^{N}$ are the lower and upper bound of $\epsilon$ at each time step. We will estimate the bounds $b_\epsilon$ of $\mathcal{W}_i$ for uncertainty quantification. To avoid an overly conservative bound of future velocity errors, we estimate the bound $b_\epsilon$ based on the past velocity errors $\epsilon_{\textrm{ini}}$ and the dynamic of the mixed traffic system.

We consider three different disturbance estimation methods: 1) zero estimation, 2) constant velocity, and 3) constant~acceleration (see Fig. \ref{fig:estimation_method} for illustration). We detail them below:  
\subsubsection{Zero estimation} 
A simplistic approach is to assume the future velocity error of the preceding vehicle is around zero, which is adopted in \method{cDeeP-LCC}~\cite{wang2023deep}. This assumption makes $\epsilon_{i,\min} \! = \! \epsilon_{i,\max} \! =\!  \mathbb{0}_{N}$ and \eqref{eqn:dDeeP} is reduced to \eqref{eqn:localDeeP} with $\epsilon_{i,\textrm{est}} \!= \! \mathbb{0}_{N}$. We consider this choice as a baseline where we only decentralize the mixed traffic system without robustification.

\subsubsection{Constant bound estimation} 
We here assume a constant velocity model, \emph{i.e.}, the velocity error of the preceding vehicle will not deviate from its current value in a short time, and the variation for the future velocity trajectory is close to its past trajectory. From the historical disturbance values $\epsilon_{i,\textrm{ini}}$, we can get the value of the current disturbance as $\epsilon_{i,\textrm{cur}} = \epsilon_{i,\textrm{ini}}(\textrm{end})$, and the disturbance variation can be estimated as $\Delta \epsilon_{i,\textrm{low}} = \min(\epsilon_{i,\textrm{ini}})-\textrm{mean}(\epsilon_{i,\textrm{ini}})$ and $\Delta \epsilon_{i,\textrm{up}} = \max(\epsilon_{i,\textrm{ini}})-\textrm{mean}(\epsilon_{i,\textrm{ini}})$. Then, the bound of the future disturbance is
    \[
    \epsilon_{i,\min} = \epsilon_{i,\textrm{cur}} + \Delta \epsilon_\textrm{low}, \  
    \epsilon_{i,\max} = \epsilon_{i,\textrm{cur}} + \Delta \epsilon_\textrm{up}.
    \]

\subsubsection{Time-varying bound estimation} 
We then assume a constant acceleration model, \emph{i.e.}, the acceleration of the preceding vehicle will not deviate significantly from its current value, and its variation in the future is close to the variation in the past. The past acceleration can be obtained from $\epsilon_{i,\textrm{ini}}$ as $a_{i,\textrm{ini}} = \frac{\epsilon_{i,\textrm{ini}}(k+1)-\epsilon_{i,\textrm{ini}}(k)}{\Delta t}$, where $\Delta t$ is the sampling time. Utilizing a similar procedure as above, we can estimate the bound of acceleration variation as $\Delta a_{i,\textrm{low}} = \min(a_{i,\textrm{ini}})-\textrm{mean}(a_{i,\textrm{ini}})$ and $\Delta a_{i,\textrm{up}} = \max(a_{i,\textrm{ini}})-\textrm{mean}(a_{i,\textrm{ini}})$. Thus, we can derive the future velocity error in time step $k$ using the bounds below:
    \[
    \begin{aligned}
    \epsilon_{i,\textrm{ini}}(\textrm{end}) + (a_{i,\textrm{cur}}+& \Delta a_{i,\textrm{low}}) \cdot k \Delta t \le \epsilon_{i}(k) \\ 
    &\le \epsilon_{i,\textrm{ini}}(\textrm{end}) + (a_{i,\textrm{cur}}+\Delta a_{i,\textrm{up}}) \cdot k \Delta t.
    \end{aligned}
    \]

\begin{figure}[t]
\centering
\includegraphics[width=0.48\textwidth]{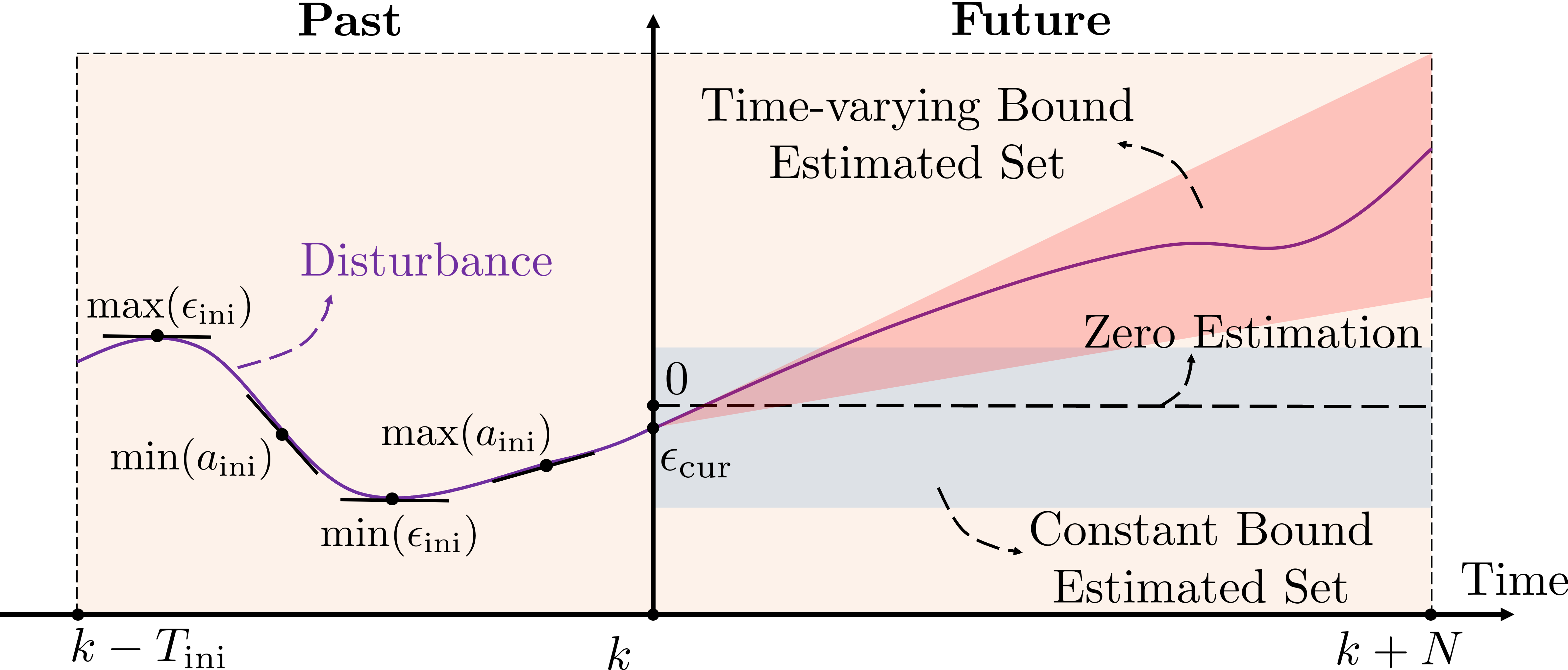}
\caption{Schematic of three disturbance estimation methods. The purple line denotes the actual disturbance trajectory, whose past is known while its future needs to be estimated. The zero estimation is denoted as the black dashed line, while the time-varying bound estimated set and the constant bound estimated set are represented as the red region and the blue region, respectively.}
\label{fig:estimation_method}
\end{figure}

Fig.~\ref{fig:estimation_method} illustrates the three different disturbance estimation methods with a sampled disturbance trajectory. In Fig.~\ref{fig:estimation_method}, we observe that a large gap can exist between the actual disturbance trajectory and the zero assumption. For the constant bound, the actual disturbance trajectory stays in the estimated region (see blue region) in the short term, but it then deviates outside the set. By contrast, the time-varying bound contains the actual trajectory in the estimated set (see red region) in this case. In most of our numerical experiments (see Section~\ref{Results}), the time-varying bounds outperform the other two methods. This is mainly because the velocity oscillations are usually low frequency but their amplitude can be high in traffic waves.  

\begin{remark}
There also exist other disturbance estimation methods in the literature~\cite{wang2022distributed,huang2021decentralized}. For distributed optimization in~\cite{wang2022distributed}, high-frequency inter-vehicle communications are required between two neighboring CF-LCC subsystems to approximate future disturbances. This might provide a more accurate estimation, but it is non-trivial to implement. For min-max \method{DeePC} in the smart grid in~\cite{huang2021decentralized}, the disturbance bound is assumed to be the largest set that contains the entire disturbance trajectory throughout the process. However, this setup is too conservative for traffic control, degrading the resulting closed-loop control performance. \hfill $\square$ 
\end{remark}

\subsection{Down-sampling Strategy}
Due to the exponential growth of the number of constraints in robust optimization~\cite{bertsimas2011theory, lofberg2012automatic}, it may not be computationally feasible to use the full estimated disturbance set \eqref{eq:uncertainty-set-W}~of dimension $N$ in problem \eqref{eqn:dDeeP}; see our derivation~of~the~constraint number in Section \ref{subsec:complexity-comp}. We here use~a~down-sampling strategy to reduce the dimension of the disturbance~set. 

The basic idea is to choose one point for every $T_\textrm{s}$ step along the $N$-dimensional disturbance trajectory and perform a linear interpolation. We denote the low-dimensional representation of the future disturbance trajectory $\epsilon_i$ as $\tilde{\epsilon}_i \in \mathbb{R}^{n_\epsilon}$ where $n_\epsilon = (\lfloor(\frac{N-2}{T_\textnormal{s}}\rfloor+2)$. The approximated representation $\hat{\epsilon}_i$ becomes:
\begin{itemize}
    \item if $1 \le k \le \tilde{k}\cdot T_\textrm{s}$:
$$
\hat{\epsilon}_i^{(k)} =
\tilde{\epsilon}_i^{(\bar{k}+1)}+((k-1) \ \textrm{mod} \ T_\textnormal{s}) \times 
\frac{\tilde{\epsilon}_i^{(\bar{k}+2)}-\tilde{\epsilon}_i^{(\bar{k}+1)}}{T_\textnormal{s}}
;
$$
    \item if $ \tilde{k}\cdot T_\textrm{s} < k \le N$:
    $$
    \hat{\epsilon}_i^{(k)} =\tilde{\epsilon}_i^{(\tilde{k}+1)}+(k-\tilde{k} \cdot T_\textrm{s}-1) \times 
\frac{\tilde{\epsilon}_i^{(\tilde{k}+2)}-\tilde{\epsilon}_i^{(\tilde{k}+1)}}{N-\tilde{k} \cdot T_\textrm{s}-1},
$$
\end{itemize}
where $\bar{k} = \lfloor\frac{k-1}{T_\textnormal{s}}\rfloor$ and $\tilde{k} = \lfloor \frac{N-2}{T_\textrm{s}} \rfloor$. Then, $\epsilon_i \in \mathbb{R}^N$ can be represented by $\tilde{\epsilon}_i \in \mathbb{R}^{n_\epsilon}$ as 
\begin{equation} \label{eqn:DownSample}
\epsilon_i \approx \hat{\epsilon}_i = E_{\epsilon} \tilde{\epsilon}_i,
\end{equation}
where $\tilde{\epsilon}_i \in \tilde{\mathcal{W}}_i$ and $\tilde{\mathcal{W}}_i$ are estimated using the same approaches as those for $\mathcal{W}_i$, and $E_{\epsilon}$ denotes the linear interpolation process above. 

As discussed in~\cite{huang2021decentralized}, using down-sampling by replacing $\epsilon_i$ with $\hat{\epsilon}_i$ may fail to cover the worst-case scenario, as the signal space of $\hat{\epsilon}_i$ is a subspace of $\epsilon_i$. Nevertheless, our extensive traffic simulations in Section~\ref{Results}  show that this approximation method could provide satisfactory traffic control performance.

\section{Tractable Reformulations and Efficient Computations}
\label{robust_computation}
Upon estimating $\tilde{W}_i$, the robust optimization problem \eqref{eqn:dDeeP} is well-defined. To run Algorithm~\ref{Alg:dDeeP-LCC}, we need to solve \eqref{eqn:dDeeP} online efficiently. In this section, we first provide a sequence of reformulations (and relaxations) for \eqref{eqn:dDeeP} since the standard solvers are not suitable to the current form. We next adapt techniques from robust optimization~\cite{lofberg2012automatic} into our context, and develop an efficient automatic transformation that reformulates each decentralized robust \method{DeeP-LCC} problem \eqref{eqn:dDeeP} into a standard conic form. With slight abuse of notations, we will use $x$ to denote the decision variable in this section, following the convention in optimization. For notational simplicity, we will also omit the subscript $i$ for the $i$-th CF-LCC subsystem when it is clear from the context.

\subsection{Reformulations via Constraint Elimination}
\label{subsec:tractable_reform}

The equality constraints \eqref{eqn:equality} in \eqref{eqn:dDeeP} make the robust optimization complicated. We first express the decision variables $g_i$ and $y_i$ in terms of $u_i$, $\sigma_{y_i}$, and $\tilde{\epsilon}_i$. Recall that $\epsilon_i$ is replaced by its low-dimensional approximation $\tilde{\epsilon}_i$ from~\eqref{eqn:DownSample}. We have 
\begin{subequations}\label{eqn:g_rep}
\begin{align}
g_i &=
H_{i}^{\dag}b \!+\!H_i^{\bot}z \label{eqn:g_rep-a},\\
y_i &= Y_{i,\textnormal{F}} g_i = Y_{i,\textnormal{F}} H_{i}^\dag b + Y_{i,\textnormal{F}}H_i^{\bot}z, \label{eqn:g_rep-b}
\end{align}
\end{subequations}
where $H_i = \textrm{col}(U_{i,\textnormal{P}}, E_{i,\textnormal{P}}, Y_{i,\textnormal{P}}, U_{i,\textnormal{F}}, E_{i,\textnormal{F}})$ denotes the Hankel data matrix for each CF-LCC system, $H_i^\dag$ is its pseudo-inverse matrix, $H^\bot = I-H^\dag H$, $z$ is an arbitrary vector in $\mathbb{R}^{T-L+1}$, and $b = \textrm{col}(u_{i,\textnormal{ini}}, \epsilon_{i,\textnormal{ini}}, y_{i,\textnormal{ini}}+\sigma_{y_i}, u_i, E_\epsilon \tilde{\epsilon}_i)$ denotes the vector on the right-hand side of \eqref{eqn:equality}. 

In the following derivation, we set $z=0$, which leads to the least-norm solution of $g_i$ for the equality constraint~\eqref{eqn:equality}. After some algebra (see Appendix~\ref{appendix:robustQP}), the min-max robust problem~\eqref{eqn:dDeeP} can be rewritten into the form of 
\begin{subequations}
\label{eqn:robustCount}
\begin{align}
\min_{\sigma_{y_i}, u_i} \ \max_{ \tilde{\epsilon}_i \in \mathcal{\tilde{W}}_i}\quad & x^\tr M x + d^\tr x + c_0 \label{eqn:Obj-robust}\\
\textrm{subject~to} \quad  
& \tilde{s}_{\min} \leq P_1 x + c_1 \leq \tilde{s}_{\max} \label{eqn:estimator-robust}, \\
& u_{\min} \le P_2 x \le u_{\max}, \label{eqn:inputlimit-robust}
\end{align}
\end{subequations}
where $x$ denotes the decision variables $\textrm{col}(u_{i}, \sigma_{y_i},  \tilde{\epsilon}_{i})$, the weight $M$ is a constant positive semidefinite matrix, $d$ is a constant vector, $c_0$ is a scalar, $P_1$, $c_1$ are constants  that construct $y_i$ in~\eqref{eqn:g_rep-b} from $x$, and $P_2$ is an index matrix extracting $u_i$ in~\eqref{eqn:inputlimit} from $x$. We note that $M$, $P_1$, $P_2$ can be pre-computed before running Algorithm~\ref{Alg:dDeeP-LCC}. Precisely, $M$, $P_1$ only depend on the Hankel matrix and the weights $\lambda_g, \lambda_y$, and $P_2$ is a constant indexing matrix. On the other hand, $d$, $c_0$ and $c_1$ depend on the initial trajectory and need to be computed at each time step. More details can be found in Appendix~\ref{appendix:robustQP}.

We then treat  $\tilde{\epsilon}_i$ as uncertainty parameters, eliminate the constant $c_0$ and transform \eqref{eqn:robustCount} into the epi-graph form below
\vspace{-5mm}
\begin{subequations}
\label{eqn:epi_robustCount}
\begin{align}
\min_{x, t} & \quad t \nonumber \\
\textrm{subject~to}& \quad  x^\tr M x + d^\tr x \leq t, \quad \forall \tilde{\epsilon}_i \in \tilde{\mathcal{W}}_i \label{eqn:obj_rb}\\ 
& \quad \tilde{s}_{\min} \le P_1 x + c_1\le \tilde{s}_{\max}, \quad \forall \tilde{\epsilon}_i \in \tilde{\mathcal{W}}_i \label{eqn:safe_rb}\\
& \quad u_{\min} \le P_2 x \le u_{\max},  %\nonumber.
\end{align}
\end{subequations}
which ensures safety constraints for all $\tilde{\epsilon}_i$ in $\tilde{\mathcal{W}}_i$.

\subsection{Vertex-based and Duality-based Strategies for Computation }
The uncertainty set $\tilde{\mathcal{W}}_i$ is a compact polytope (see Section \ref{sec:dist_est}), and thus \eqref{eqn:obj_rb} and \eqref{eqn:safe_rb} have an infinite number of constraints, which is not directly solvable in the current form. We here adapt techniques from robust optimization~\cite{lofberg2012automatic} to solve \eqref{eqn:epi_robustCount} by exploiting the representations of $\tilde{\mathcal{W}}_i$. 

In particular, it is known that the compact polytope $\tilde{\mathcal{W}}_i$ has two different representations: 1) internal representations as its convex hull of extreme points (\emph{i.e.}, its vertices), and 2) external representations as an intersection of affine subspaces. We proceed to derive a computable form for each of them, leading to two different methods respectively: 1) the vertex-based strategy and 2) the duality-based strategy. 

To simplify notations, we denote 
\begin{subequations}
    \begin{align}
        \tilde{\mathcal{W}}_i &= \mathrm{conv} (w_1, \ldots, w_{n_\textrm{v}})  \label{eq:extreme-point-Wi}\\
        \tilde{\mathcal{W}}_i &= \{\tilde{\epsilon}_i \in \mathbb{R}^{n_\epsilon} \mid \tilde{A}_{\epsilon} \tilde{\epsilon}_i \leq \tilde{b}_{\epsilon}\},  
        \label{eq:affine-subspace-Wi}
    \end{align}
\end{subequations}
where $n_\textrm{v}$ is the number of extreme points, $w_1, \ldots, w_{n_\textrm{v}} \in  \mathbb{R}^{n_\epsilon}$ are the extreme points of $\tilde{\mathcal{W}}_i$, and $\tilde{A}_\epsilon = \textrm{col}(I_{n_\epsilon},-I_{n_\epsilon})$, $\tilde{b}_{\epsilon} = \textrm{col}(\tilde{\epsilon}_{i,\max}, \tilde{\epsilon}_{i,\min})$ are the low-dimensional approximation for affine constraints in \eqref{eq:uncertainty-set-W}.

The first method is summarized in the following proposition. 
\begin{proposition}[Method I: Vertex-based Strategy] 
\label{prop:Vetex}Using the extreme point representation \eqref{eq:extreme-point-Wi}, the optimization \eqref{eqn:epi_robustCount} is equivalent to the following problem
\begin{subequations} \label{finalformVertex}
\begin{align}
\min_{x, t} \quad& t \nonumber\\
\textrm{subject to} \quad & x_j^\tr M x_j + d^\tr x_j \le t,   j = 1,\ldots, n_{\textrm{v}}, \label{FFVertexC1}\\ 
& \tilde{s}_{\min} \le P_1 x_j +c_1 \le \Tilde{s}_{\max}, j = 1,\ldots, n_{\textrm{v}}, \label{FFVertexC2} \\
& u_{\min} \le P_2 x \le u_{\max} \label{FFVertexC3},
\end{align}
\end{subequations}
where $x_j$ denotes the decision variables $\textrm{col}(u_{i}, \sigma_{y_i},  w_j)$ by fixing $\tilde{\epsilon}_i$ to be one of the extreme points $w_j$.  
\end{proposition}

\begin{IEEEproof}
    We only need to prove that \eqref{eqn:obj_rb} and \eqref{eqn:safe_rb} are equivalent to \eqref{FFVertexC1} and \eqref{FFVertexC2}, respectively. 

The equivalence \eqref{eqn:safe_rb} $\Leftrightarrow$ \eqref{FFVertexC2} is a direct consequence of the basic result in linear programs (LPs). In particular, the right-hand inequality of \eqref{eqn:safe_rb} is the same as 
\begin{equation}
\label{safe_opt}
\begin{aligned}
 \tilde{s}_\textnormal{max} \geq \max_{\tilde{\epsilon}_i \in \tilde{\mathcal{W}}_i} \;\; p_l ^ \tr x  + c_{1,l} , \qquad l = 1, \ldots, N,
\end{aligned}
\end{equation}
where $p_l^\tr$ is the $l$-th row vector in $P_1$, and $c_{1,l}$ is $l$-th element in $c_1$. Recall that $N$ is the predication horizon, and \eqref{safe_opt} just enforces the spacing error of the $i$-th CAV to be upper bounded by the maximum spacing error $\tilde{s}_\textnormal{max}$ for all time steps. Since~\eqref{safe_opt} is an LP with a bounded feasible region, the maximum is always attained at one of its vertices of $\tilde{\mathcal{W}}_i$, represented in \eqref{eq:extreme-point-Wi}. Therefore, \eqref{safe_opt} is equivalent to the right-hand inequality of \eqref{FFVertexC2}. Similarly, the argument for the left-hand inequality of \eqref{eqn:safe_rb} can be obtained. 

The equivalence \eqref{eqn:obj_rb} $\Leftrightarrow$ \eqref{FFVertexC1} requires a slightly different argument from convex sets. The direction \eqref{eqn:obj_rb} $\Rightarrow$ \eqref{FFVertexC1} is immediate since  all vertices are contained in set $\tilde{\mathcal{W}}_i$. As for the converse direction  \eqref{FFVertexC1} $\Rightarrow$ \eqref{eqn:obj_rb},  we consider any element $\tilde{\epsilon}_i \in \tilde{\mathcal{W}}_i$. A basic result in convex analysis guarantees the convex representation 
$$
\tilde{\epsilon}_i = \sum_{j=1}^{n_{\textrm{v}}} \alpha_j w_j, \quad \sum_{j=1}^{n_\textrm{v}} \alpha_j = 1, \quad \alpha_j \geq 0,
$$ 
where $w_j, j = 1, \ldots ,n_{\textrm{v}}$ are the extreme points \eqref{eq:extreme-point-Wi}.  
Following the same notation in \eqref{finalformVertex}, we also have $x = \sum_{j=1}^{n_{\textrm{v}}} \alpha_j x_j$. Then, we have 
$$
\begin{aligned}
    x^\tr M x + d^\tr x  &\leq \sum_{j=1}^{n_{\textrm{v}}} \alpha_j ( x_j^\tr M x + d^\tr x_j ) \\
    & \leq \sum_{j=1}^{n_{\textrm{v}}} \alpha_j t  = t, \quad \forall \tilde{\epsilon}_i \in \tilde{\mathcal{W}}_i,
\end{aligned}
$$
where the first inequality comes from the fact that $x^\tr M x + d^\tr x$ is a convex function, and the second inequality is from~\eqref{FFVertexC1}.  This completes the proof.
\end{IEEEproof}

Problem~\eqref{finalformVertex} is a standard convex problem and can be solved using existing solvers. However, the number of extreme points $n_{\textrm{v}}$ can be large even for a simple polytope, which indicates the number of constraints in~\eqref{finalformVertex} can be large.~We next introduce another equivalent formulation that utilizes~duality analysis to reduce the number of constraints, which is summarized below. 

\begin{proposition}[Method II: Duality-based strategy]
    Using the representation as an intersection of affine subspaces~\eqref{eq:affine-subspace-Wi}, the optimization~\eqref{eqn:epi_robustCount} is equivalent to the following problem
    \begin{subequations} \label{finalformDual}
\begin{align}
\min_{x_{\textrm{d}}, t, \lambda_1, \lambda_2} \quad& t \nonumber\\
\textrm{subject to} \quad %& x_j^\tr M x_j + d^\tr x_j \le t,  j = 1,\ldots, n_{\textrm{v}}, \label{FFDualC1}\\ 
& p_{l, \textrm{d}}^{\tr} x_{\textrm{d}}  + \tilde{b}_\epsilon^\tr \lambda_{l,1} +  c_{1,l}\le 
 \tilde{s}_{\max}, \label{FFDualC1}\\
 & \tilde{A}_\epsilon^\tr \lambda_{l,1} - p_{l,\epsilon} = 0,\label{FFDualC2} \\ 
&-p_{l, \textrm{d}}^\tr x_{\textrm{d}}  + \tilde{b}_\epsilon^\tr \lambda_{l,2} -  c_{1,l} \le -\tilde{s}_{\min},\label{FFDualC3} \\
 &\tilde{A}_\epsilon^\tr \lambda_{l,2} + p_{l,\epsilon} = 0,\label{FFDualC4} \\
&\lambda_{l,1} \geq 0, \lambda_{l,2} \geq 0, \ l = 1,2,\ldots, N,  \label{FFDualC5} \\
& \eqref{FFVertexC1}, \eqref {FFVertexC3}, \nonumber 
\end{align}
\end{subequations}
where $x_{\textrm{d}}=\textrm{col}(u_i, \sigma_{y_i})$ denotes the decision variable, and $\lambda_{l,1} \in \mathbb{R}^{2n_{\epsilon}}$ is the dual variable corresponding to the affine constraints in \eqref{safe_opt}, with $\lambda_1 = \textrm{col}(\lambda_{1,1} \lambda_{2,1}, \ldots, \lambda_{N,1})$. Similar definitions hold for $\lambda_{l,2}$ and $\lambda_2$. The coefficients $\tilde{A}_\epsilon, \tilde{b}_\epsilon$ can be found in~\eqref{eq:affine-subspace-Wi}, and $c_{1,l}, p_{l}$ are the same as~\eqref{safe_opt}. We use  $p_{l,\textrm{d}}$ to denote $\textrm{col}(p_{l,u}, p_{l,\sigma_y})$ with $p_l$ being partitioned as $\textrm{col}(p_{l,u}, p_{l,\sigma_{y}},p_{l,\epsilon})$, corresponding to the original variables $u_i, \sigma_{y_i}$ and $\tilde{\epsilon}_i$, respectively. 
\end{proposition}

\begin{IEEEproof}
We have shown that~\eqref{eqn:obj_rb} is equivalent to~\eqref{FFVertexC1}. We only need to prove the equivalence between \eqref{eqn:safe_rb} and \eqref{FFDualC1} - \eqref{FFDualC5}. Recall that the right-hand inequality of \eqref{eqn:safe_rb} is equivalent to the LP in~\eqref{safe_opt}. 
Via duality analysis, we shall prove that \eqref{safe_opt} will be the same as 
\eqref{FFDualC1}, \eqref{FFDualC2} and \eqref{FFDualC5}. The argument for the left-hand inequality of \eqref{eqn:safe_rb} is identical. 

To derive the dual of the LPs in \eqref{safe_opt}, we first write its Lagrangian function for $l = 1,\ldots, N$ as 
\[
\begin{aligned}
L(\tilde{\epsilon}_i, \lambda_{l, 1}) &= p_l^\tr x + c_{1,l} + \lambda_{l,1}^\tr (\tilde{b}_\epsilon - \tilde{A}_\epsilon \tilde{\epsilon}_i) \\
&=(p_{l,\epsilon}^\tr-\lambda_{l,1}^\tr \tilde{A}_\epsilon)\tilde{\epsilon}_i + p_{l,\textrm{d}}^\tr x_\textrm{d} + \lambda_{l,1}^\tr \tilde{b}_\epsilon + c_{1,l}.\\  
\end{aligned}
\]
The dual function is thus given by
$$
h(\lambda_{l,1}) = \begin{cases} p_{l,\textrm{d}}^\tr x_\textrm{d} + \lambda_{l,1}^\tr \tilde{b}_\epsilon + c_{1,l} & \text{if} \;\; p_{l,\epsilon}^\tr-\lambda_{l,1}^\tr \tilde{A}_\epsilon = 0\\
-\infty & \text{otherwise}.    
\end{cases}
$$
Then, the dual problem for the LPs in $\eqref{safe_opt}$ becomes 
\begin{equation}
\label{eqn:safe_dual}
\begin{aligned}
f_l^\star = \min_{\lambda_{l,1}} \quad & p_{l, \textrm{d}}^{\tr} x_{\textrm{d}} + \tilde{b}_\epsilon^\tr \lambda_{l,1} + c_{1,l} \\
\textrm{subject~to} \quad &  \tilde{A}_\epsilon^\tr \lambda_{l,1}  - p_{l,\epsilon} = 0, \lambda_{l,1} \succeq 0. 
\end{aligned}
\end{equation}
The strong duality of LPs ensures that \eqref{safe_opt} is the same as 
\begin{equation} \label{eq:dual-bound}
 \tilde{s}_\textnormal{max} \geq f_l^\star, \qquad l = 1, \ldots, N.
\end{equation}

We can thus replace the right-hand inequality of~\eqref{eqn:safe_rb} with~\eqref{eqn:safe_dual} and~\eqref{eq:dual-bound}, leading to a bi-level minimization problem.~Simple algebra allows us to combine both levels\footnote{This operation is standard. The interested reader may see, \emph{e.g.}, Section 2.1 of this note: \url{https://zhengy09.github.io/ECE285/lectures/L17.pdf}.}, confirming the equivalence between~\eqref{eqn:safe_dual} and~\eqref{FFDualC1},~\eqref{FFDualC2},~\eqref{FFDualC5}. 
This completes our proof.
\end{IEEEproof}

\subsection{Complexity Comparison and Automatic Transformation}
\label{subsec:complexity-comp}

We here compare the complexity of Methods I and II. One main complexity difference is caused by using different representations of \eqref{safe_opt} to derive~\eqref{eqn:safe_rb}, which corresponds to \eqref{FFVertexC2} in Method I and \eqref{FFDualC1}-\eqref{FFDualC5} in Method II. 

In our problem, the uncertainty set is a polytope with dimension $n_\epsilon$ and we have $n_\textrm{v} = 2^{n_\epsilon}$. For Method~I, \eqref{FFVertexC2} represents $2N \cdot 2^{n_\epsilon}$ inequality constraints and its complexity will increase exponentially as $n_\epsilon$ grows up. Moreover, $2N$ will be a relatively large coefficient since the prediction horizon in our problem is usually larger than $20 $ (corresponding to $1\,\mathrm{s}$). For Method II,~\eqref{FFDualC1} - \eqref{FFDualC5} together represent $2N\cdot(3n_\epsilon+1)$ inequality constraints, whose size is much smaller than the $2N\cdot 2^{n_\epsilon}$ constraints in~\eqref{FFVertexC2} given a large value of $n_\epsilon$. Table~\ref{Table:complexity} illustrates the complexity of~\eqref{finalformVertex} and~\eqref{finalformDual}. As we can observe, when the prediction horizon $N$ is fixed, Method II has a smaller complexity than Method I given a large value of $n_\epsilon$.

\begin{table}[t]
\caption{Complexity comparison between \eqref{finalformVertex} and \eqref{finalformDual}.}
\begin{threeparttable}
\setlength{\tabcolsep}{4pt}
    \begin{tabular}{ccc}
         \toprule
         & Decision Variables Number& Constraints Number \\
         \midrule
         \textbf{Method I}  & $(m_i+2)T_\textrm{ini}+N+1$ & $2^{n_\epsilon}+N \cdot 2^{n_\epsilon+1}+2N$\\
         \textbf{Method II} & $(m_i+2)T_\textrm{ini}+N+1 + 4 N n_\epsilon$ & $2^{n_\epsilon} + 2N(3n_\epsilon + 2)$\\
         \bottomrule
    \end{tabular}
    \begin{tablenotes}
    \footnotesize
    \item[1] We recall that $m_i$ is the number of HDVs in the $i$-th subsystem, $T_\textrm{ini}$ is the length of the initial trajectory, $N$ is the length of prediction horizon and $n_\epsilon$ is the dimension of the approximated disturbance set. 
    \end{tablenotes}
    \label{Table:complexity}
    \end{threeparttable}
\end{table}

The theoretical complexity result is only one key factor that influences the time consumption. For numerical computation, we need to convert \eqref{finalformVertex} or \eqref{finalformDual} into a standard conic optimization before passing them to a numerical solver.~For this process, we find that modeling packages such as YALMIP~\cite{Lofberg2004} often introduce too much overhead time consumption. Also, similar instances of \eqref{finalformVertex} or \eqref{finalformDual} need to be converted repeatedly for the receding horizon online control in Algorithm \ref{Alg:dDeeP-LCC}. Thus, we here 
implement an automatic transformation from \eqref{finalformVertex} or \eqref{finalformDual} into standard conic optimization, so that they can be directly solved by solvers (\emph{e.g.}, Sedumi~\cite{Strum1999}, Mosek~\cite{mosek}). 

Precisely, their standard conic dual form for \eqref{finalformVertex} or \eqref{finalformDual} can be written as 
\begin{equation}
\label{eqn:Standard_Conic}
    \begin{aligned}
        \max_{y,s} \quad & b^\tr y \\
     \mathrm{subject~to} \quad  & A^\tr y + s = c \\
        &s \in \mathcal{K},
    \end{aligned}
\end{equation}
where $A$ is a matrix, $y,b,c,s$ are vectors and $\mathcal{K}$ is a cone. We present a detailed transformation from \eqref{finalformVertex} and \eqref{finalformDual} to \eqref{eqn:Standard_Conic} in Appendix \ref{appendix:conic_form}. We note that our automatic transformation does not introduce extra decision variables or constraints.

\begin{figure}[t]
\centering
\vspace{-3mm}
\setlength{\abovecaptionskip}{3pt}
\subfigure[Modeling time]{\includegraphics[width=0.24\textwidth]{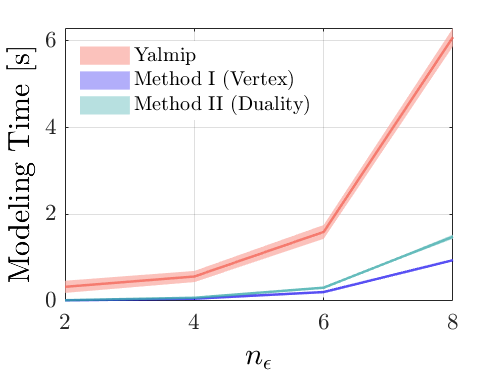}}
\subfigure[Solver time]{\includegraphics[width=0.24\textwidth]{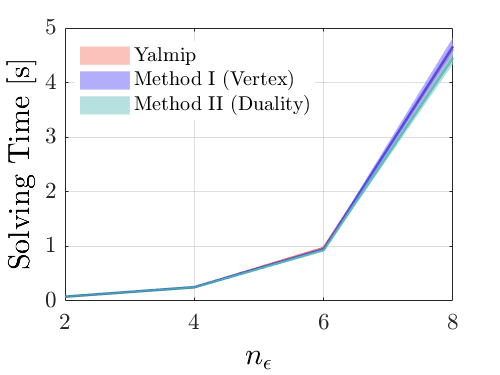}}
\caption{Comparison of computation time for different methods.}
\label{fig:computation_time_comp}
\vspace{-2mm}
\end{figure}

Here, we use a simple numerical experiment to demonstrate the advantage of our automatic transformation compared to modeling packages such as YALMIP~\cite{Lofberg2004} for solving \eqref{finalformVertex} and \eqref{finalformDual}. In particular, we consider two types of time costs: modeling time and solver time. Modeling time refers to the time needed to update the parameters of the conic form \eqref{eqn:Standard_Conic} in each iteration, while solver time is the time the solver takes to solve the resulting conic program \eqref{eqn:Standard_Conic}. 
For the numerical experiment, we fix the parameters $m_i = 3, T_\textrm{ini} = 20, N = 50$ and vary $n_\epsilon$. For each value of $n_\epsilon$, we compute $20$ instances using different methods. The results are shown in Fig. \ref{fig:computation_time_comp}.~As we can observe, both Vertex-based and Duality-based methods have comparable computational performance. With the increase of $n_\epsilon$, Method I has faster modeling time while Method II performs slightly better in the solver time. Also, our automatic transformation provides much faster modeling time with similar solving time compared to YALMIP. In our following experiments (Section \ref{Results}), we will thus use the automatic transformation for Method I to solve \method{dDeeP-LCC}~\eqref{eqn:dDeeP}.

%\balance
\section{Numerical Experiments} \label{Results}
In this section, we present three nonlinear and non-deterministic traffic simulations to validate the performance of \method{dDeeP-LCC}\footnote{All our implementation and experimental scripts can be found at \url{https://github.com/soc-ucsd/Decentralized-DeeP-LCC}.}. The driving behaviors of HDVs are modeled by a nonlinear OVM model, similar to~\cite{wang2021controllability}. A noise with a uniform distribution of $\mathbb{U}[-0.1, 0.1]$ $\mathrm{m/s}^2$ is added to the acceleration signal of each HDV in the simulation. We utilize \method{dDeeP-LCC}(Zero), \method{dDeeP-LCC}(Constant), and \method{dDeeP-LCC}(Time-varying) to represent different disturbance estimation methods in Section \ref{sec:dist_est}.

% For the mixed traffic system in Fig.~\ref{fig:MixTrafSys}, 
We consider a system with sixteen vehicles behind one head vehicle, among which there are four CAVs and twelve HDVs (\emph{i.e.}, $n=16$, $q=4$). The CAVs are located in the third, sixth, tenth, and thirteenth vehicles, respectively (\emph{i.e.}, $\mathcal{S} = \{3, 6, 10, 13\}$). Accordingly, there are four CF-LCC subsystems, and the corresponding index sets for HDVs are $F_1 = \{4, 5\}$, $F_2 = \{7, 8, 9\}$, $F_3 = \{11, 12\}$ and $F_4 = \{14, 15, 16\}$. Other parameters are set as follows:
\begin{itemize}
    \item \textit{Offline data collection:} lengths of pre-collected trajectories are chosen as $T=700$ for a small data set and $T = 1500$ for a large data set with $\Delta t = 0.05$ $\mathrm{s}$. They are collected around the equilibrium state of the system with velocity 15 $\mathrm{m/s}$. We use a uniform distributed signal of $\mathbb{U}[-1,1]$ to generate both $u^\D$ and $\epsilon^\D$ which satisfies persistent excitation requirement in Propositions \ref{proposition:data-rep-central} and \ref{proposition:data-rep-decen}.
    \item \textit{Online predictive control:} The prediction horizon and the initial signal sequence are set to $N = 50$ and $T_\textnormal{ini} = 20$, respectively. In the cost function \eqref{eqn:objFun}, we have $\lambda_{g_i} = 10$, $\lambda_{y_i} = 10000$, $R_i = 0.1I_N$, $Q_i = I_N \otimes \textrm{diag}(Q_{v_i}, w_{s_i})$ where $Q_{v_i} = \textrm{diag}(1,\ldots, 1) \in \mathbb{R}^{m_i+1}$ and $w_{s_i} = 0.5$. The spacing constraints for the CAVs are set as $s_{\max} = 40$ $\mathrm{m}$, $s_{\min} = 5$ $\mathrm{m}$ and the bound of the spacing error is updated at each time step as $\tilde{s}_{\max} = s_{\max} - s^*$ and $\tilde{s}_{\min} = s_{\min} - s^*$. The acceleration bound of the CAVs is set to $a_{\max} = 2$ $\mathrm{m/s}^2$ and $a_{\min} = -5$ $\mathrm{m/s}^2$. Note that $s^*$ is updated in each iteration of Algorithm \ref{Alg:dDeeP-LCC} based on the current equilibrium state, estimated by the past trajectory of the head vehicle (as suggested in~\cite{wang2023deep}).
\end{itemize}

For comparison with model-based methods, we include centralized MPC and decentralized MPC which share the same parameter settings corresponding to \method{cDeeP-LCC} and \method{dDeeP-LCC} (\emph{e.g.}, $R_i, Q_i, \tilde{s}_{\max}, \tilde{s}_{\min}, a_{\max}, a_{\min}$) except the predictive model. The centralized MPC utilizes the parametric model \eqref{eqn:ModelMixTraffic} for prediction and we denote it as \method{cMPC}. The decentralized MPC for each CF-LCC subsystem using the parametric model \eqref{eqn:ModelCFLCC} for prediction and we further integrate it with the different disturbance estimation approaches which are represented as \method{dMPC}(Zero), \method{dMPC}(Constant) and \method{dMPC}(Time-varying). We note that both \eqref{eqn:ModelMixTraffic} and \eqref{eqn:ModelCFLCC} are generally unknown while we use the actual model in the simulation.

\begin{figure}[t]
\centering
\setlength{\abovecaptionskip}{2pt}
\subfigure[All HDVs]{\includegraphics[width=0.235\textwidth]{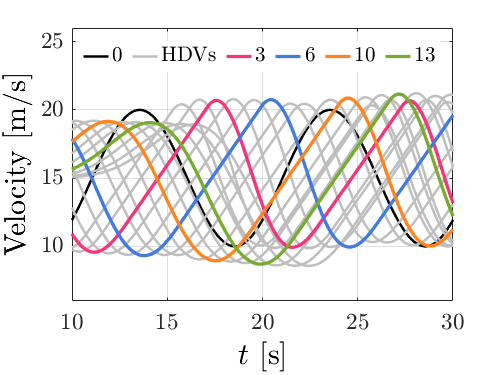} \label{subfig:Sin_HDV}} \vspace{-2mm}\\
\subfigure[\method{cDeeP-LCC}]{\includegraphics[width=0.235\textwidth]{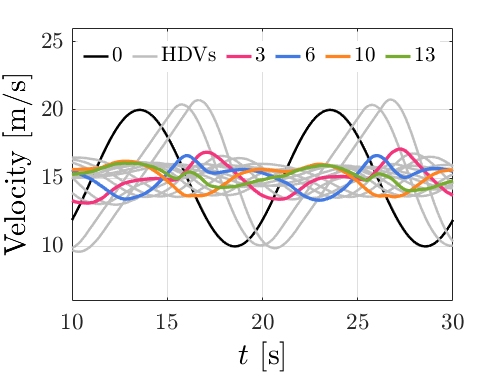} \label{subfig:Sin_cDeeP}} 
\subfigure[\method{dDeeP-LCC}]{\includegraphics[width=0.235\textwidth]{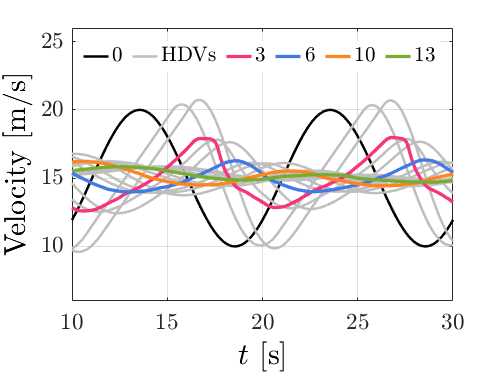}\label{subfig:Sin_dDeeP}} 

\caption{Velocity profiles in Experiment A where the head vehicle is under sinusoidal perturbation. The black profile represents the head vehicle and the gray profile represents the HDVs. The red profile (vehicle 3), blue profile (vehicle 6), orange profile (vehicle 10) and green profile (vehicle 13) represent the first to the fourth CAV, respectively. (a) All vehicles are HDVs. (b) CAVs utilize \method{cDeeP-LCC}. (c) CAVs utilize \method{dDeeP-LCC}.}
\label{fig:Sin_Vel}
\end{figure}

\subsection{Performance Validation around an Equilibrium State}

Our first experiment (Experiment A) simulates a traffic wave scenario by imposing a sinusoidal perturbation on the head vehicle. This is to validate the performance of \method{dDeeP-LCC} in smoothing traffic waves around a fixed equilibrium point. The head vehicle accelerates and decelerates periodically around the equilibrium velocity $15\,\mathrm{m/s}$ and its velocity profile is $v_{\textrm{head}}(t) = (15+5 \cdot \sin(0.2\pi t))$ $\mathrm{m/s}$ (see black profile in Fig.~\ref{fig:Sin_Vel}). To quantify the closed-loop performance, we consider 
the mean squared velocity error (MSVE), defined as
\[
\text{MSVE} = \frac{\Delta t}{n(t_f-t_0)}\sum_{t = t_0}^{t_f}\sum_{i=1}^{n}(v_i(t)-v_0(t))^2,
\]
where $t_0$ and $t_f$ denote the start and end time of the simulation. 

As shown in Fig. \ref{subfig:Sin_HDV}, when all vehicles are HDVs, the amplitude of velocity perturbation from the head vehicle is amplified along the vehicle string. In contrast, the amplitude of these perturbations is greatly reduced when the CAVs are equipped with \method{cDeeP-LCC} and it decreases progressively along the vehicle string for CAVs utilizing \method{dDeeP-LCC}(Time-varying), as shown in Fig. \ref{subfig:Sin_cDeeP} and Fig.~\ref{subfig:Sin_dDeeP}, respectively. Also, \method{dDeeP-LCC}(Time-varying) reduces $91.8\%$ of the MSVE, which is comparable to $93.8\%$ for \method{cDeeP-LCC}; see Table~\ref{table:msve}. This indicates that \method{dDeeP-LCC} with an appropriate disturbance estimation method can effectively dissipate traffic waves with negligible performance degradation compared with \method{cDeeP-LCC}. The minor improved performance of \method{cDeeP-LCC} comes from it takes the entire mixed traffic system into consideration and does not need to estimate the coupling dynamics between CF-LCC subsystems. Thanks to the decentralization, one major benefit of \method{dDeeP-LCC} is its computational scalability (see details in Section~\ref{subsection:time-consumption}).

\begin{table}[t]
\caption{MSVE Statistics in Experiment A}
\centering
\begin{threeparttable}
\setlength{\tabcolsep}{2pt}
    \begin{tabular}{cccccc}
    \toprule
      \multirow{2}*{}& \multirow{2}*{HDVs} & \multirow{2}*{\method{cDeeP-LCC}} &\multicolumn{3}{c}{\method{dDeeP-LCC}} \\
     \noalign{\vskip -2pt}
     \cmidrule(lr{0.5em}){4-6}
     \noalign{\vskip -2pt}
     & & & Zero & Constant & Time-varying \\
    \hline
    MSVE & 9.25 & 0.57 & 1.28 & 1.33 & 0.76 \\
    \hline
    Reduction& N/A & \textbf{93.8\%} & 86.2\%& 85.6\%& \textbf{91.8\%}\\
    \noalign{\vskip -2pt}
    \bottomrule
    \end{tabular}
    \begin{tablenotes}
    \footnotesize
    \item[1] Units for MSVE is $\textrm{m}^2/\textrm{s}^2$ in this table. 
    \end{tablenotes}
    \end{threeparttable}
    \label{table:msve}
\end{table}

\subsection{Traffic Improvement in Comprehensive Simulations}
\label{subsec:compreh-sim}

\begin{figure*}[t]
\centering
\setlength{\abovecaptionskip}{0pt}
\subfigure[\method{cDeeP-LCC}]{\includegraphics[width=0.49\textwidth]{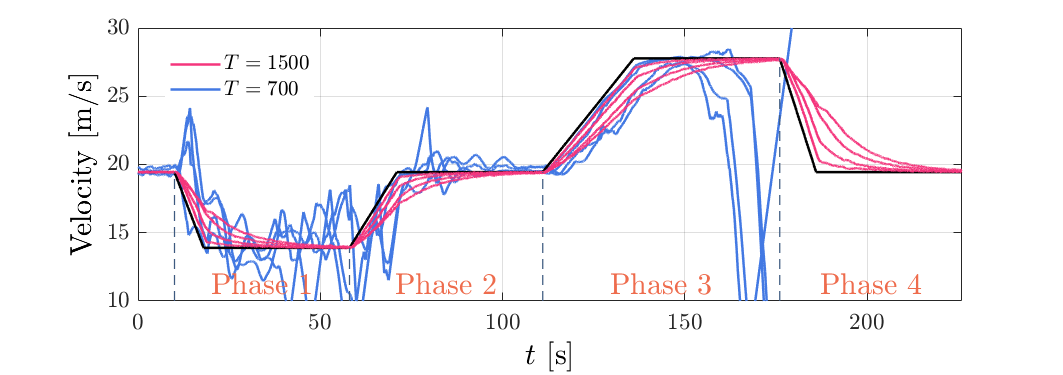}\label{subfig:NEDC_Cen}} 
\subfigure[\method{dDeeP-LCC}(Zero)]{\includegraphics[width=0.49\textwidth]{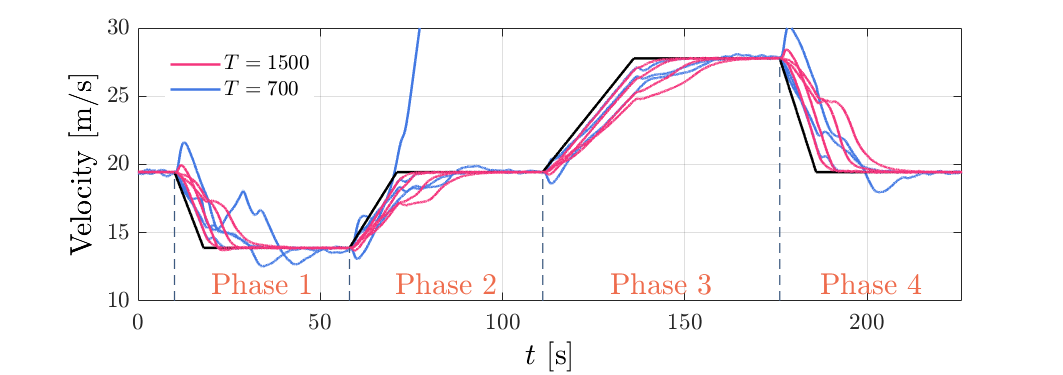}\label{subfig:NEDC_zero}}
\subfigure[\method{dDeeP-LCC}(Bound)]{\includegraphics[width=0.49\textwidth]{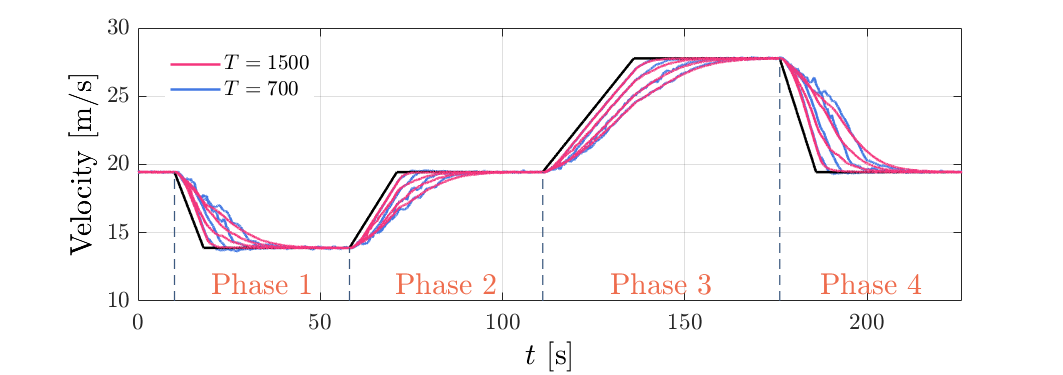}\label{subfig:NEDC_constant}} 
\subfigure[\method{dDeeP-LCC}(Time-varying)]{\includegraphics[width=0.49\textwidth]{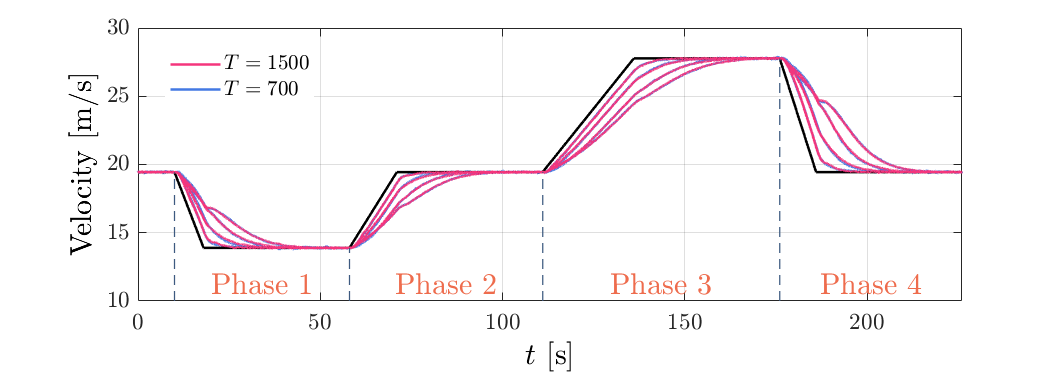}\label{subfig:NEDC_timev}}
\caption{Velocity profiles in Experiment B. The black profile denotes the head vehicle. The red profiles and the blue profiles represent the data-driven control for CAVs with data sets of size $T=1500$ and $T=700$, respectively. (a) The CAV utilizes \method{cDeeP-LCC}. (b)(c)(d) The CAV utilizes \method{dDeeP-LCC} with zero, constant bound and time-varying bound estimation approaches.}
\label{fig:nedc_vel}
\end{figure*}

Our next Experiment B validates the performance of \method{dDeeP-LCC} in a comprehensive simulation scenario with time-varying equilibrium states. In particular, we consider the New European Driving Cycle (NEDC)~\cite{DieselNet2013}, where the velocity profile for the head vehicle is shown in the black profile in Fig. \ref{fig:nedc_vel}.  We also compute the fuel consumption of $4$ CF-LCC subsystems for comparison. The fuel consumption rate $f_i$ ($\mathrm{mL/s}$) for the $i$-th vehicle is calculated by~\cite{bowyer1985guide}
\[
f_i = \left\{
\begin{aligned}
    & 0.444+0.090 R_i v_i + [0.054a_i^2v_i]_{a_i>0}, & \text{if} \ R_i >0,\\
    & 0.444, & \text{if} \ R_i \le 0,
\end{aligned}
\right.
\]
where $a_i$ represents the acceleration of vehicle $i$ and $R_i = 0.333+0.00108v_i^2+1.200a_i$; see~\cite{bowyer1985guide} for details. 

Fig. \ref{fig:nedc_vel} shows the velocity profiles for \method{cDeeP-LCC} and \method{dDeeP-LCC} utilizing different estimation methods with two data sets ($T=700, 1500$). When using the large data set ($T=1500$), all the methods enable the CAVs to track the desired velocity and smooth traffic flow (red curves in Fig. \ref{fig:nedc_vel}). By contrast, when using the small data set ($T=700$),  both \method{cDeeP-LCC} and \method{dDeeP-LCC}(Zero) fail to stabilize the closed-loop traffic (blue curves in Fig. \ref{subfig:NEDC_Cen} \ref{subfig:NEDC_zero}). Still, \method{dDeeP-LCC}(Constant) and \method{dDeeP-LCC}(Time-varying) can stabilize the mixed traffic system. Moreover, \method{dDeeP-LCC}(Time-varying) achieves the best car-following behaviors for the CAVs with the smallest velocity oscillations (see the blue curves in Fig.~\ref{subfig:NEDC_timev}). Although \method{dDeeP-LCC} does not require communications between different decentralized units, the consensus of the mixed traffic system is achieved (see the convergence behavior in Fig. \ref{fig:nedc_vel}) since all CF-LCC subsystems aim to track the same velocity. Additionally, consensus can also be achieved for a given reference velocity, or by making each CF-LCC subsystem track the velocity of its preceding vehicle, thanks to the cascading structure.

Compared with \method{dDeeP-LCC}(Time-varying), the performance degradation of the other methods is caused by the behavior representation error due to a smaller data set and the disturbance estimation error. Both types of errors can lead to a mismatch between the online prediction and real system behavior. We compute the theoretical minimum required data length via \eqref{eqn:data-length}, which is $611$ for a centralized setting and $233 $ for a decentralized setting (computed using the largest subsystem).
Since we use a linearized model to approximate the nonlinear system, the data length needs to be much larger than the theoretical value. Thus, \method{cDeeP-LCC} relies on a relatively large data set ($T=1500$) for a valid representation of the entire system behavior, while the small data set ($T=700$) is not large enough which fails to stabilize the closed-loop performance. In \method{dDeeP-LCC}, from the zero bound, the constant bound to the time-varying bound, the estimation error decreases and provides more margin for potential representation errors, leading to improved control performance. Thus, \method{dDeeP-LCC}(Time-varying) outperforms other methods for a small~data~set. 

\begin{table*}[t]
\caption{Fuel Consumption Comparison (Unit: $\mathrm{mL}$)}
\centering
\begin{threeparttable}
\setlength{\tabcolsep}{3pt}
    \begin{tabular}{ccccccccccccc}
    \toprule
    \noalign{\vskip -2pt}
     \multirow{2}*{}& \multirow{2}*{HDVs} & \multirow{2}*{\method{cDeeP-LCC}} &\multicolumn{3}{c}{\method{dDeeP-LCC}} &\multicolumn{3}{c}{\method{dDeeP-LCC} (0.2\,$\mathrm{s}$ delay)} & \multirow{2}*{\method{cMPC}} & \multicolumn{3}{c}{\method{dMPC}} \\
     \noalign{\vskip -2pt}
     \cmidrule(lr{0.5em}){4-6}\cmidrule(lr{0.5em}){7-9} \cmidrule(lr{0.5em}){11-13}
     \noalign{\vskip -2pt}
     & & & Zero & Constant & Time-varying & Zero & Constant & Time-varying&  & Zero & Constant & Time-varying \\
     \noalign{\vskip -2pt}
    \midrule
     NEDC & 6548.36 & 6346.16 & 6388.47 & 6360.69 & 6356.14 & 6387.67 & 6360.74 & 6356.57 & 6344.58 & 6343.42 & 6366.39 & 6358.50 \\
    \midrule
     Braking & 1296.15 & 894.68 & 900.36 & 880.35 & 874.35 & 909.02 & 884.25 & 876.74 & 880.11 & 900.16 & 884.94 & 869.49 \\
    \bottomrule
    \end{tabular}
    \end{threeparttable}
    \label{table:msve_fuel_res}
\end{table*}

The fuel consumption results with the large data set ($T=1500$) are listed in column $2$ to $8$ of Table~\ref{table:msve_fuel_res}, and a detailed comparison between \method{cDeeP-LCC} and \method{dDeeP-LCC}(Time-varying) is shown in Table \ref{Table:fuel_comp}. Compared with the case with all HDVs, all control methods can reduce fuel consumption, and a certain degree of communication delay (\emph{i.e.}, 0.2\,$\mathrm{s}$) has little effect on the performance of \method{dDeeP-LCC} as the fuel consumption remains similar with the delayed signal. Also, we observe that \method{dDeeP-LCC}(Time-varying) performs the best among different disturbance estimation methods and the improvement of implementing \method{DeeP-LCC} in the braking phase (Phases 1 and 4) is higher than the accelerating phases (Phases 2 and 3) as shown in Table \ref{Table:fuel_comp}. Moreover, for the overall reduction of fuel consumption, the \method{dDeeP-LCC}(Time-varying) achieves comparable fuel economy as the \method{cDeeP-LCC}, which is $2.94\%$ and $3.09\%$. Although the \method{cDeeP-LCC} performs slightly better than \method{dDeeP-LCC} because of more accurate modeling of dynamic couplings, it requires more offline collected data as well as longer computation time.  

\begin{remark}[Comparison with model-based methods]
    The fuel consumption results for \method{cMPC} and \method{dMPC} are also shown in Table \ref{table:msve_fuel_res}. We first observe that \method{cDeeP-LCC} has comparable performance as \method{cMPC} and the maximum difference of fuel consumption between \method{dDeeP-LCC}(Time-varying) and \method{dMPC} is $0.66\%$ among all cases. Moreover, for both data-driven and model-based predictive controllers, utilizing time-varying disturbance estimation provides superior performance when the velocity of the head vehicle changes rapidly. It is worth noting that \method{dDeeP-LCC} does not need prior knowledge of the system model, and it directly synthesizes the optimal controller \eqref{eqn:dDeeP} by utilizing a single trajectory while model-based controllers require to conduct system identification. 
\end{remark}
\begin{remark}[Robustification for centralized controllers]
 Both \method{cMPC} and \method{cDeeP-LCC} can be integrated with different disturbance estimations. However, the disturbance estimation method has a negligible effect for centralized controllers, and the differences of fuel consumption among various estimation approaches are smaller than $0.13\%$ as shown in Table \ref{table:fuel-central}. 
 \begin{table}[t!]
\caption{Fuel Consumption Comparison for Centralized Controllers(Unit: $\mathrm{mL}$)}
\centering
\begin{threeparttable}
\setlength{\tabcolsep}{2pt}
    \begin{tabular}{ccccccc}
    \toprule
    \noalign{\vskip -2pt}
     \multirow{2}*{}& \multicolumn{3}{c}{\method{cMPC}} &\multicolumn{3}{c}{\method{cDeeP-LCC}} \\
     \noalign{\vskip -2pt}
      \cmidrule(lr{0.5em}){2-4} \cmidrule(lr{0.5em}){5-7}
     \noalign{\vskip -2pt}
     &Zero & Constant &Time-varying &Zero & Constant & Time-varying \\
     \noalign{\vskip -2pt}
    \midrule
     NEDC & 6344.58 & 6344.28 & 6344.07 & 6346.16 & 6345.90 & 6346.10\\
     Braking & 880.11  & 878.99 & 878.97 & 894.68 & 895.36 & 894.88\\
    \bottomrule
    \end{tabular}
    \end{threeparttable}
    \label{table:fuel-central}
    \vspace{-3pt}
\end{table}
 The main reason is the coupling dynamics between subsystems (\emph{i.e.}, velocity error of preceding vehicles of subsystems) are system states in the centralized setting and have been predicted via the system model and considered by the optimal controller. Moreover, we shall illustrate that the computation time of the \method{cDeeP-LCC} with robustification can increase significantly in Section~\ref{subsection:time-consumption}. Since the robustification of the centralized controller costs much more computational resources with little performance increase, we omit their results in Table~\ref{table:msve_fuel_res}.
\end{remark}
\begin{table}[t]
\caption{Fuel Consumption Comparison in Experiment B (unit: $\mathrm{mL}$)}
\vspace{-2mm}
\begin{threeparttable}
    \begin{tabular}{c c c c}
         \toprule
         \noalign{\vskip-2pt}
          &All HDVs& \method{cDeeP-LCC}& \method{dDeeP-LCC} \\
          \noalign{\vskip-2pt}
         \midrule
         Phase 1 & 749.23  &  700.99 ($\downarrow 6.44\%$) & 699.89 ($\downarrow 6.58\%$)\\
         Phase 2 & 1351.18  & 1331.75 ($\downarrow 1.44\%$)  & 1333.29 ($\downarrow 1.32\%$)\\
         Phase 3 & 2904.48   & 2882.21 ($\downarrow 0.77\%$)  & 2886.86 ($\downarrow 0.61\%$)\\
         Phase 4 & 1304.59   & 1191.91 ($\downarrow 8.64\%$)  & 1196.34 ($\downarrow 8.30\%$)\\
        Total Process & 6548.36  &6346.16 ($\downarrow$ \textbf{3.09\%})   & 6356.14 ($\downarrow$ \textbf{2.94\%})\\
         \bottomrule
    \end{tabular}
    \begin{tablenotes}
    \footnotesize
    \item[1] We use time-varying bound estimation for \method{dDeeP-LCC} in this table. 
    \end{tablenotes}
    \end{threeparttable}
    \label{Table:fuel_comp}
    \vspace{-6mm}
\end{table}

\subsection{Safety Performance in Emergence Braking}
\label{subsec:safety}
Our last Experiment C uses an emergency braking scenario to validate the safety performance of \method{dDeeP-LCC} with different disturbance estimation methods. The head vehicle that moves at $15\,\mathrm{m}/\mathrm{s}$ will brake with maximum deceleration $-5\,\mathrm{m}/\mathrm{s}^2$, remain at $5 \,\mathrm{m/s}$ for a while and accelerate back to $15 \,\mathrm{m}/\mathrm{s}$. We carry out the same experiment for all the \method{dDeeP-LCC} controllers with $100 $ small data sets ($T=700$) and $100$ large data sets ($T=1500$). As discussed in Section~\ref{subsec:uncert-quant}, we take the zero estimation as a baseline since it can be considered as applying \method{cDeeP-LCC} for each subsystem without robustification. We recall that the safety spacing constraint for the CAV is set from $5\,\mathrm{m}$ to $40\,\mathrm{m}$. The ``violation" and ``emergency" are defined as cases where at least one CAV's spacing deviates more than $1 \,\mathrm{m}$ and $5\,\mathrm{m}$ from the safe range, respectively. We note that one of the following three undesired situations will happen when an emergency occurs: 1) A rear-end collision happens; 2) The CAV's spacing is too large which decreases the traffic capacity; 3) The controller can not stabilize the mixed traffic~system. 

The safety performance results are listed in Table \ref{Table:safe_rate}. With large data sets, \method{dDeeP-LCC}(Zero) has a high violation rate and emergency rate, which are $88\%$ and $68\%$ respectively, while the other two methods can provide much better safety guarantee ($0\%$ violation rate). The zero estimation method shows the worst robustness performance against large disturbances, which leads to a high risk of rear-end collisions. Given the small data sets, the violation rate and emergency rate for both \method{dDeeP-LCC}(Constant) and \method{dDeeP-LCC}(Time-varying) remain $0\%$, but they are close to $100\%$ for \method{dDeeP-LCC}(Zero). Moreover, \method{dDeeP-LCC}(Time-varying) performs the best in reducing the fuel consumption among all data-driven methods with the large data set ($T=1500$) as shown in Table \ref{table:msve_fuel_res} which decreases $32.5\%$ fuel cost compared with the case of all HDVs. These results validate the superior performance of \method{dDeeP-LCC}(Time-varying), which works for small data sets and can also provide the best safety performance. 
\begin{table}[t]
\centering
\caption{Collision and safety Constraint Violation Rate}
\vspace{-2mm}
\begin{threeparttable}
\setlength{\tabcolsep}{2pt}
    \begin{tabular}{c c c c c c c}
         \toprule
         \noalign{\vskip-2pt}
          & \multicolumn{2}{c}{Zero}&\multicolumn{2}{c}{Constant}&\multicolumn{2}{c}{Time-varying} \\
          \cmidrule(lr{0.5em}){2-3} \cmidrule(lr{0.5em}){4-5} \cmidrule(lr{0.5em}){6-7}
          & $T\!= \!700$ & $T\!= \!1500$ & $T\!= \!700$ & $T\!= \!1500$ & $T\!= \! 700$ & $T\!= \!1500$\\
          \noalign{\vskip-2pt}
         \midrule
         Violation & 99$\%$ & 88$\%$ & $\mathbf{0}\%$ & 0$\%$ & $\mathbf{0}\%$ & 0$\%$\\ 
         Emergency & 97$\%$ & 68$\%$ & $\mathbf{0}\%$ & 0$\%$ & $\mathbf{0}\%$ & 0$\%$\\
         \bottomrule
    \end{tabular}
    \end{threeparttable}
    \label{Table:safe_rate}
    \vspace{-2mm}
\end{table}

\begin{figure}[t]
\centering
\subfigure[Large offline data set with $T=1500$]{\includegraphics[width=0.24\textwidth]{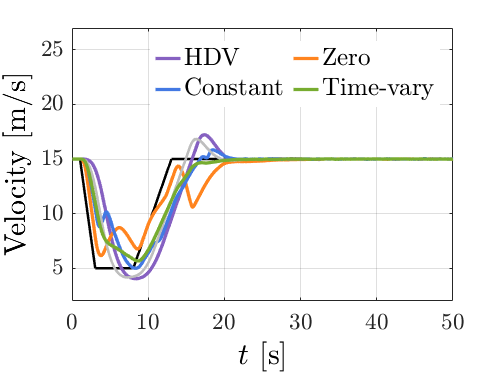}
\includegraphics[width=0.24\textwidth]{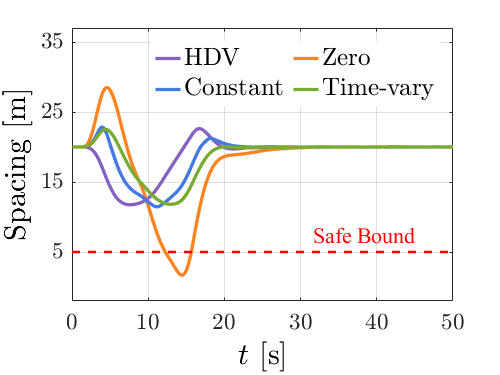}\label{fig:space_T=1500}}

\subfigure[Small offline data set with $T=700$]{\includegraphics[width=0.24\textwidth]{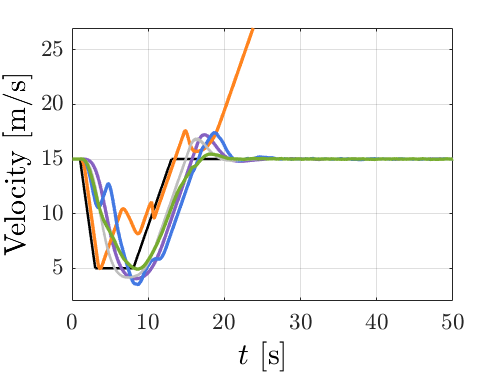}
\includegraphics[width=0.24\textwidth]{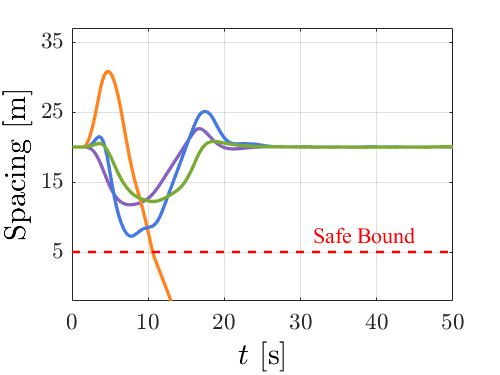}\label{fig:space_T=700}}
\vspace{-3.5mm}
\caption{Simulation results in Experiment C. The black profile and the gray profile denote the head vehicle and the preceding vehicle, respectively. The orange profile, blue profile and the green profile represent \method{dDeeP-LCC}(Zero), \method{dDeeP-LCC}(Constant) and \method{dDeeP-LCC}(Time-varying) respectively, while the purple profile corresponds to the all HDV case. (a) and (b) show the velocity and spacing profiles at different sizes of data sets. }
\label{fig:safety_compare}
\vspace{-2.5mm}
\end{figure}

In Fig.~\ref{fig:safety_compare}, we further illustrate the trajectories of the first CAV  to compare the performance of utilizing different disturbance estimation methods. Compared with the case of all human drivers, we observe that all methods using a large data set can smooth the traffic flow. However, using a small data set, \method{dDeeP-LCC} with both zero and constant bound estimation approaches leads to undesired velocity oscillations. For the spacing profiles, \method{dDeeP-LCC}(Time-varying) and \method{dDeeP-LCC}(Constant) satisfy the safety constraint for both large and small data sets, but \method{dDeeP-LCC}(Zero) is likely to cause a rear-end collision at a small data set. Although the safety constraint is imposed in \method{dDeeP-LCC}(Zero), the mismatch between the prediction and real system behavior leads to its failure in the implementation. This is due to the oversimplified assumption in  \method{dDeeP-LCC}(Zero) that the preceding vehicle tracks the equilibrium velocity accurately. On the other hand, the \method{dDeeP-LCC}(Time-varying) relies on a more accurate estimation of the future disturbances, and introduces a robust design that requires safety satisfaction for all the possible future disturbances. Therefore, \method{dDeeP-LCC}(Time-varying) provides a much better safety guarantee with a smoother velocity profile.

\subsection{Computational Efficiency of \method{dDeeP-LCC}} \label{subsection:time-consumption}
Finally, we compare computation time for each method. In our experiments, we solve~\eqref{eqn:cDeeP-LCC} and~\eqref{finalformVertex} using Mosek~\cite{mosek} on an Intel Core i7-9750H CPU with 16GB RAM. We consider~\eqref{eqn:cDeeP-LCC} as a specific form of~\eqref{finalformVertex} by setting $n_\textrm{v} = 1$ and $w_1 = \mathbb{0}_{n_\epsilon}$, and utilize the same method to solve them for comparison. We use \method{cDeeP-LCC}(Zero), \method{cDeeP-LCC}(Constant) and \method{cDeeP-LCC}(Time-varying) to represent solving the optimization problem in the centralized setting with different disturbance estimation approaches. % with no advanced solving techniques implemented.

As shown in Table~\ref{table:computation_time_Res}, the $\method{dDeeP-LCC}$ with three disturbance estimation methods achieve much faster computation time compared with $\method{cDeeP-LCC}$. The proposed \method{dDeeP-LCC}(Time-varying) method is $85.1\%$ faster in terms of average computation time compared with the \method{cDeeP-LCC}(Zero) and demonstrates superior control and safety performance as previously discussed. Note that all the CAVs can solve their control inputs individually in the decentralized setting, and thus the computational cost is related to the size of the subsystem and will not increase as the size of the entire system grows. As the number of subsystems increases, the computational cost for $\method{cDeeP-LCC}$ will keep increasing, while it will remain nearly unchanged for $\method{dDeeP-LCC}$. Moreover, the large dimension of the entire mixed traffic system can significantly increase the complexity of the optimization problem \eqref{eqn:dDeeP} after robustification. The computation time for \method{cDeeP-LCC} after robustification is longer than $1.5\,\mathrm{s}$ (see columns 3 and 4 in Table~\ref{table:computation_time_Res}) which is unacceptable for real implementation. This reveals the scalability of the decentralized approach in terms of both computation and data privacy, which supports real-time implementation for large-scale mixed traffic systems.

\begin{table*}[t]
\centering
\caption{Computation Time (seconds)}
\label{table:computation_time_Res}
\vspace{-2.5mm}
\begin{threeparttable}
    \begin{tabular}{cccccccc}
    \toprule
    \noalign{\vskip -2pt}
     &\multicolumn{4}{c}{\method{cDeeP-LCC}} &\multicolumn{3}{c}{\method{dDeeP-LCC}} \\
     \noalign{\vskip -2pt}
     \cmidrule(lr{0.5em}){2-5} \cmidrule(lr{0.5em}){6-8}
     \noalign{\vskip -2pt}
     & ADMM & Zero & Constant & Time-varying &Zero & Constant & Time-varying \\
     \noalign{\vskip -2pt}
    \midrule
    \noalign{\vskip -2pt}
    Traffic Wave & 0.0848 &0.41 & 1.69& 1.69 & 0.030 & 0.059 &  0.059\\
    NEDC & 0.0895&0.37 & 1.55 & 1.50 & 0.032& 0.059& 0.057\\
    Braking & 0.0732 &0.39 &1.66 &1.64 & 0.027 &0.057 &0.059  \\
    \noalign{\vskip -2pt}
    \bottomrule
    \end{tabular}
    \end{threeparttable}
    \vspace{-2mm}
\end{table*}

We further implement the controller proposed in \cite{wang2022distributed} to conduct the computation performance comparison with using ADMM to accelerate \method{cDeeP-LCC} which is denoted as \method{cDeeP-LCC}(ADMM). As shown in Table \ref{table:computation_time_Res}, the computation time of \method{dDeeP-LCC} is also shorter than \method{cDeeP-LCC}(ADMM). In the parallel computation of the \method{cDeeP-LCC}(ADMM), its subproblem is constructed based on the information of each CF-LCC subsystem so that the complexity of the subproblem is similar as \method{dDeeP-LCC}. Moreover, the ADMM algorithm requires impractical frequent communication between subsystems.

\section{Conclusions} \label{Conclusion}
In this paper, we have proposed a decentralized robust data-driven predictive control for CAVs, called \method{dDeeP-LCC}, to smooth mixed traffic flow. We first decouple a mixed~traffic system into several cascading CF-LCC subsystems, where the dynamical coupling is treated as a bounded disturbance. We formulate $\method{dDeeP-LCC}$ as a robust min-max optimization problem against the worst disturbance, where three estimation methods are discussed for the disturbance bound. We have also developed an efficient computational method to solve the \method{dDeeP-LCC} for online predictive control. Traffic simulations confirm that \method{dDeeP-LCC} achieves comparable control performance with centralized \method{DeeP-LCC} in~\cite{wang2023deep}, while providing better safety guarantees, achieving shorter computation time, and naturally preserving data privacy. 

Interesting future directions include 1) adapting different regularization strategies and data-preprocessing methods from~\cite{dorfler2022bridging,shang2023convex,zhang2023dimension} for further improved performance, 2)~developing direct data-driven methods for nonlinear traffic control (such as~\cite{berberich2022linear}) which bypass the local linearization around the equilibrium state and improve the model accuracy, and 3) considering the dynamically changing structure of the mixed traffic system which requires the online data recollection process. Practical experimental validations similar to~\cite{wang2022implementation} and more large-scale numerical experiments are also interesting future topics.

\bibliographystyle{IEEEtran}
\bibliography{reference}

\appendix

In this appendix, we first present the transformation of our decentralized robust \method{DeeP-LCC} \eqref{eqn:dDeeP} into the robust quadratic optimization problem \eqref{eqn:robustCount}. We then provide standard conic forms of the vertex-based strategy (Method I) \eqref{finalformVertex} and duality-based strategy (Method II) \eqref{finalformDual} for directly solving the decentralized robust \method{DeeP-LCC} using solvers, without utilizing time-consuming modeling packages. We further illustrate the essential of ensuring the persistent excitation condition (see Definition \ref{def:PE}) and discuss the implementation of \method{dDeeP-LCC}.

% \section{Calculation details}
\subsection{Quadratic reformulation} \label{appendix:robustQP}
We show the process of getting $Q$, $P_1$, $P_2$, $d$, $c_0$ and $c_1$ in~\eqref{eqn:robustCount}. The main idea is substituting \eqref{eqn:g_rep} into \eqref{eqn:dDeeP} and using $x$ to replace variables $u_i, \sigma_{y_i}$ and $\tilde{\epsilon}_i$. 

We first derive the objective function of \eqref{eqn:Obj-robust}. After substituting, the objective function of \eqref{eqn:dDeeP} is
\begin{equation}
\label{eqn:quad_obj}
\begin{aligned}
&||u_i||^2_{R_i} + ||y_i||^2_{Q_i} +  \lambda_{g_i} ||g_i||_2^2 + \lambda_{y_i} ||\sigma_{y_i}||_2^2 \\ 
= \; & u_i^\tr R_i u_i + b^\tr S_1 b + \lambda_{g_i}b^\tr S_2 b + \lambda_{y_i} \sigma_{y_i}^\tr \sigma_{y_i}.
\end{aligned}
\end{equation}
where $S_1 = H_i^{\dag \tr} Y_{i,\textrm{F}}^\tr Q_i Y_{i,\textrm{F}} H_i^{\dag}$ and $S_2 =  H_i^{\dag \tr} H_i^{\dag}$ are symmetric positive semidefinite matrix. We denote $u_i, \sigma_{y_i}$~and~$b$~as 
\[
\begin{aligned}
& u_i = F_1 x, \quad \sigma_{y_i} = F_2 x, \quad  b = b_{\textrm{ini}} + F_3 x,
\end{aligned}
\]
where $b_{\textrm{ini}} = \textrm{col}(u_{i,\textrm{ini}}, \epsilon_{i,\textrm{ini}}, y_{i,\textrm{ini}}, \mathbb{0}, \mathbb{0})$ and
\[
\begin{aligned}
& F_1 = \begin{bmatrix}
    I_{N} & \mathbb{0} & \mathbb{0}
\end{bmatrix}, \qquad  F_2 = \begin{bmatrix}
    \mathbb{0} & I_{(m_i+2)T_{\textrm{ini}}} & \mathbb{0}
\end{bmatrix}, \ 
\\
& F_3 = \begin{bmatrix}
    \mathbb{0} & \mathbb{0} & \mathbb{0} \\
    \mathbb{0} & \mathbb{0} & \mathbb{0} \\
    \mathbb{0} & I_{(m_i+2)T_{\textrm{ini}}} & \mathbb{0} \\
    I_{N} & \mathbb{0} & \mathbb{0} \\
    \mathbb{0} & \mathbb{0} & E_\epsilon
\end{bmatrix}.
\end{aligned}
\] Then \eqref{eqn:quad_obj} could be derived as 
\[
\begin{aligned}
&u_i^\tr R_i u_i + b^\tr S_1 b + \lambda_{g_i}b^\tr S_2 b + \lambda_{y_i} \sigma_{y_i}^\tr \sigma_{y_i} \\
=\;&  x^\tr F_1^\tr R_i F_1 x + (b_{\textrm{ini}} + F_3 x)^\tr S_1 (b_{\textrm{ini}} + F_3 x) \\
&\quad + \lambda_{g_i}(b_{\textrm{ini}} + F_3 x)^\tr S_2 (b_{\textrm{ini}} + F_3 x) + \lambda_{y_i} x^\tr F_2^\tr F_2^\tr x \\
=\;& x^T(F_1^\tr R_i F_1 + F_3^\tr S_1 F_3 + \lambda_{g_i}F_3^\tr S_2 F_3 + \lambda_{y_i} F_2^\tr F_2^\tr)x \\
& \quad + 2(b_{\textrm{ini}}^\tr S_1 F_3 + \lambda_{g_i}b_{\textrm{ini}}^\tr S_2 F_3)x + b_{\textrm{ini}}^\tr (S_1+\lambda_{g_i}S_2) b_{\textrm{ini}},
\end{aligned}
\]
so that we can get
\[
\begin{aligned}
M &= F_1^\tr R_i F_1 + F_3^\tr S_1 F_3 + \lambda_{g_i}F_3^\tr S_2 F_3 + \lambda_{y_i} F_2^\tr F_2^\tr,\\
d &= 2(F_3^\tr S_1 b_{\textrm{ini}}+ \lambda_{g_i}F_3^\tr S_2 b_{\textrm{ini}}),\\
c_0 &=  b_{\textrm{ini}}^\tr (S_1+\lambda_{g_i}S_2) b_{\textrm{ini}}.
\end{aligned}
\]

We then show the derivation for constraints. For the safety constraint \eqref{eqn:estimator-robust}, after substituting \eqref{eqn:g_rep-b} into \eqref{eqn:safety} and using the same decomposition of $b$ in the previous part, we have
\[
\begin{aligned}
&\tilde{s}_{\min} \le G_1 Y_{i,\textrm{F}} H_i^\dag(b_{\textrm{ini}} + F_3 x) \le \tilde{s}_{\max} \\ 
\Leftrightarrow \ & \tilde{s}_{\min} \le G_1 Y_{i,\textrm{F}} H_i^\dag F_3 x + G_1 Y_{i,\textrm{F}} H_i^\dag b_{\textrm{ini}} \le \tilde{s}_{\max}.
\end{aligned}
\]
Thus, we have
\[
P_1 = G_1 Y_{i,\textrm{F}} H_i^\dag F_3, \quad c_1 = G_1 Y_{i,\textrm{F}} H_i^\dag b_{\textrm{ini}}.
\]
For the input limitation \eqref{eqn:inputlimit-robust}, we need to extract $u_i$ from $x$ and we have $P_2 = F_1$. 

\subsection{Conic form transformation}
\label{appendix:conic_form}
We show how to change~\eqref{finalformVertex} and~\eqref{finalformDual} to their conic form and get corresponding $A, b, c, \mathcal{K}$. The key point is to derive each constraint to its corresponding conic inequality. We use subscript ``$\textrm{s}", ``\textrm{o}", ``\textrm{z}$" to represent parameters corresponding to the second-order cone, non-negative orthant and zero cone. Recall that $x_j$ denotes fixing $\tilde{\epsilon}_i$ in $x$ with one of extreme points~$w_j$ and $x_\textrm{d}$ is defined as $\textrm{col}(u_i, \sigma_{y_i})$. 

\subsubsection{Method I}
For Method I, the decision variable $y$ is  $\textrm{col}(t, u_i, \sigma_{y_i})$ and it is obvious that $b = \textrm{col}(1, \mathbb{0}, \mathbb{0})$ for the objective function. There are two types of conic inequalities at the constraints which are second-order cone and non-negative~orthant. 

We first derive constraint~\eqref{FFVertexC1} to its second order cone form. It is normal to transform the convex quadratic constraints to the following form 
\begin{equation}
\label{eqn:rotated_second_order_cone}
\begin{aligned}
    \norm*{ \begin{bmatrix}
        2\Gamma x_j \\
        t-d^T x_j-1
    \end{bmatrix}}_2\le t-d^T x_j+1, \ j = 1,2,\ldots n_\textrm{v},
\end{aligned}
\end{equation}
where $\Gamma$ is gotten from decomposing $M$ as $\Gamma^\tr \Gamma$ using the fact that $M$ is a positive semidefinite matrix. We further divide $\Gamma$ and $d$ as $\Gamma = \begin{bmatrix} \Gamma_\textrm{d}, \Gamma_w\end{bmatrix}$ and $d = \textrm{col}(d_\textrm{d}, d_w)$ which corresponds to the multiplication with $x_d$ and $w_j$. Then \eqref{eqn:rotated_second_order_cone} can be derived as the second order cone form 
\[
    \begin{bmatrix}
    t-d_\textrm{d}^\tr x_d -d_w^\tr w_j +1 \\
        2(\Gamma_d x_d + \Gamma_w w_j) \\
        t-d_\textrm{d}^\tr x_d -d_w^\tr w_j -1
    \end{bmatrix} \in \mathcal{K}_{\textrm{s}j}, \quad j = 1,\ldots, n_\textrm{v},
\]
and the corresponding $A_{\textrm{s}j}$ and $c_{\textrm{s}j}$ are 
\[
A_{\textrm{s}j}^{\tr} = \begin{bmatrix}
  -1 & d_\textrm{d}^\tr\\
   \mathbb{0} & -2 \Gamma_\textrm{d} \\ 
  -1 & d_\textrm{d}^\tr
\end{bmatrix}, \
c_{\textrm{s}j} = \begin{bmatrix}
    -d_w^\tr w_j+1 \\
    2 \Gamma_w w_j \\ 
    -d_w^\tr w_j -1
\end{bmatrix}, \ j= 1,\ldots, n_\textrm{v}.
\]

We then derive constraints~\eqref{FFVertexC2} and~\eqref{FFVertexC3} to their non-negative orthant form. For constraint~\eqref{FFVertexC2}, we could represent $P_1$ as $[P_{1\textrm{d}}, P_{1 w}]$ and each evaluation of \eqref{FFVertexC2} in $w_j$ can be derived as 
\[
\begin{bmatrix}
    s_{\max} -P_{1\textrm{d}} x_\textrm{d} -  P_{1 w} w_j -c_1  \\
    -s_{\min} +P_{1\textrm{d}} x_\textrm{d} +  P_{1 w} w_j +c_1
\end{bmatrix} \in \mathcal{K}_{\textrm{o}j},
\]
and the corresponding $A_{\textrm{o}j}$ and $c_{\textrm{o}j}$ for $j = 1,\ldots, n_\textrm{v}$ are
\[
\begin{aligned}
A_{\textrm{o}j}^{\tr} = \begin{bmatrix}
    \mathbb{0} & P_{1\textrm{d}} \\
    \mathbb{0} & -P_{1\textrm{d}}
\end{bmatrix}, \quad
c_{\textrm{o}j} = \begin{bmatrix}
    s_{\max} - P_{1w} w_j-c_1 \\
    -s_{\min} + P_{1w} w_j + c_1
\end{bmatrix}.
\end{aligned}
\]
For the constraint~\eqref{FFVertexC3}, similar to derivation of \eqref{FFVertexC2}, we have
$
\begin{bmatrix}
    u_{\max} - u_i \\
    u_i - u_{\min}
\end{bmatrix}  \in  \mathcal{K}_{\textrm{o}u}
$
and the corresponding $A_{\textrm{o}u}$ and $c_{\textrm{o}u}$ are
\[
A_{\textrm{o}u}^{\tr} = \begin{bmatrix}
    \mathbb{0} & I & \mathbb{0} \\
    \mathbb{0} & -I & \mathbb{0}
\end{bmatrix}, \quad
c_{\textrm{o}u} = \begin{bmatrix}
    u_{\max} \\
    -u_{\min}
\end{bmatrix}.
\]

The final form of $A$ and $c$ for \eqref{finalformVertex} are gotten from stacking matrices derived above and we have
\[
\begin{aligned}
A & = \textrm{col}(A_{\textrm{o}u}, A_{\textrm{o}1},\ldots, A_{\textrm{o}n_\textrm{v}},A_{\textrm{s}1},\ldots,A_{\textrm{s}n_\textrm{v}}),\\
c & = \textrm{col}(c_{\textrm{o}u}, c_{\textrm{o}1},\ldots, c_{\textrm{o}n_\textrm{v}},c_{\textrm{s}1},\ldots,c_{\textrm{s}n_\textrm{v}}).
\end{aligned}
\]
Cone $\mathcal{K}$ is the product of $n_{\textrm{v}}+1$ non-negative orthants and $n_{\textrm{v}}$ second-order cone which is 
\[
\mathcal{K} = \mathcal{K}_{\textrm{o}u} \times \mathcal{K}_{\textrm{o}1} \times \cdots \mathcal{K}_{\textrm{o}n_{\textrm{v}}} \times \mathcal{K}_{\textrm{s}1} \times \cdots \times \mathcal{K}_{\textrm{s}n_{\textrm{v}}}.
\]
\subsubsection{Method II}
For Method II, the decision variable is $y = \textrm{col}(t_i, u_i, \sigma_{y_i}, \lambda_1, \lambda_2)$ and $b$ in the objective function is $\textrm{col}(1, \mathbb{0}, \mathbb{0}, \mathbb{0}, \mathbb{0})$. There are three types of conic inequalities at the constraints which are second-order cone, non-negative orthant and zero cone. 

We omit the detailed derivation for the same constraints \eqref{FFVertexC1} and \eqref{FFVertexC2} in \eqref{finalformDual}. We only need to add an extra zero block column at the right side of each $A_{\textrm{s}j}$ and $A_{\textrm{o}u}$ because of additional decision variables $\lambda_1$ and $\lambda_2$. Forms of $A_{\textrm{s}j}$ for $j = 1,\ldots, n_\textrm{v}$ and $A_{\textrm{o}u}$ can be represented as
\[
A_{\textrm{s}j}^{\tr} = \begin{bmatrix}
  -1 & d_\textrm{d}^\tr & \mathbb{0}\\
  \mathbb{0} & -2 \Gamma_\textrm{d} & \mathbb{0}\\ 
  -1 & d_\textrm{d}^\tr & \mathbb{0}
\end{bmatrix}, \quad A_{\textrm{o}u}^{\tr} = \begin{bmatrix}
    \mathbb{0} & I & \mathbb{0} & \mathbb{0}\\
    \mathbb{0} & -I & \mathbb{0} & \mathbb{0}
\end{bmatrix}.
\]
We will have the same $\mathcal{K}_{\textrm{s}j}, c_{\textrm{s}j}, j = 1,\ldots, n_{\textrm{v}}$ and $\mathcal{K}_{\textrm{s}u}, c_{\textrm{o}u}$ as in Method I.  

We then present the conic derivation for $\eqref{FFDualC1}$-$\eqref{FFDualC5}$. For \eqref{FFDualC1}, \eqref{FFDualC3}, \eqref{FFDualC5}, they could be represented by non-negative orthants which are
\[
\begin{bmatrix}
 \tilde{s}_{\max} -p_{l, \textrm{d}}^{\tr} x_{\textrm{d}}  - \tilde{b}_\epsilon^\tr \lambda_{l,1} -  c_{1,l} \\
 -\tilde{s}_{\min} +p_{l, \textrm{d}}^\tr x_{\textrm{d}}  - \tilde{b}_\epsilon^\tr \lambda_{l,2} +  c_{1,l} 
\end{bmatrix} \in \mathcal{K}_{\textrm{o}l}, \quad l = 1,\ldots, N,
\] and
$
\begin{bmatrix}
  \lambda_1 \\
  \lambda_2
\end{bmatrix} \in \mathcal{K}_{\textrm{o}\lambda}.
$
The corresponding $A_{\textrm{o}l}$ and $c_{\textrm{o}l}$ for $\mathcal{K}_{\textrm{o}l}$ are
\[
\begin{aligned}
A_{\textrm{o}l}^\tr &= \begin{bmatrix}
    0 & p_{l,d}^\tr & \mathbb{0}_{1 \times 2n_\epsilon(l-1)} & \tilde{b}_\epsilon^\tr & \mathbb{0} \\
    0 & -p_{l,d}^\tr & \mathbb{0}_{1 \times 2n_\epsilon(N+l-1)} & \tilde{b}_\epsilon^\tr & \mathbb{0}
\end{bmatrix}, \\
c_{\textrm{o}l}^\tr & = \begin{bmatrix}
    \tilde{s}_{\max} - c_{1,l} \\
     -\tilde{s}_{\min} - c_{1,l}
\end{bmatrix}, \quad l = 1,\ldots, N,
\end{aligned}
\]
and the corresponding $A_{\textrm{o}\lambda}$ and $c_{\textrm{o}\lambda}$ for $\mathcal{K}_{\textrm{o}\lambda}$ are
\[
A_{\textrm{o}\lambda}^\tr = \begin{bmatrix}
    \mathbb{0} & -I & \mathbb{0}
\end{bmatrix}, \quad c_{\textrm{o}\lambda} = \mathbb{0}.
\]
For \eqref{FFDualC2} and \eqref{FFDualC4}, they can be considered as zero cones which are 
\[
\begin{bmatrix}
\tilde{A}_\epsilon^\tr \lambda_{l,1} - p_{l,\epsilon} \\ 
\tilde{A}_\epsilon^\tr \lambda_{l,2} + p_{l,\epsilon}
\end{bmatrix} \in \mathcal{K}_{\textrm{z}l}, \quad l = 1,\ldots, N,
\]
and $A_{\textrm{z}l}$ and $c_{\textrm{z}l}$ correspond to them are
\[
\begin{aligned}
&A_{\textrm{z}l}^\tr = \begin{bmatrix}
    \mathbb{0} & \mathbb{0}_{n_\epsilon \times 2n_\epsilon(l-1)} & \tilde{A}_\epsilon^\tr & \mathbb{0} \\
    \mathbb{0} & \mathbb{0}_{n_\epsilon \times 2n_\epsilon(N+l-1)} & \tilde{A}_\epsilon^\tr & \mathbb{0}
\end{bmatrix}, \quad c_{\textrm{z}l} = \begin{bmatrix}
    p_{l,\epsilon} \\
    -p_{l, \epsilon}
\end{bmatrix}.
\end{aligned}
\]

The final form of $A$ and $c$ for \eqref{finalformDual} could be represented as 
\[
\begin{aligned}
A = \textrm{col}(& A_{\textrm{z}1},\ldots,A_{\textrm{z}N}, A_{\textrm{o}\lambda}, A_{\textrm{o}u}, \\
&A_{\textrm{o}1},\ldots, A_{\textrm{o}N}, A_{\textrm{s}1},\ldots, A_{\textrm{s}n_\textrm{v}}), \\
c = \textrm{col}(& c_{\textrm{z}1},\ldots,c_{\textrm{z}N}, c_{\textrm{o}\lambda}, c_{\textrm{o}u}, \\
&c_{\textrm{o}1},\ldots, c_{\textrm{o}N}, c_{\textrm{s}1},\ldots, c_{\textrm{s}n_\textrm{v}}). \\
\end{aligned}
\]
Cone $\mathcal{K}$ is the product of $N$ zero cones, $N+2$ non-negative orthants and $n_v$ second-order cones, which is 
\[
\begin{aligned}
\mathcal{K} = & \ \mathcal{K}_{\textrm{z}1}\times \cdots \times \mathcal{K}_{\textrm{z}N} \times \mathcal{K}_{\textrm{o}\lambda} \times \mathcal{K}_{\textrm{o}u} \times \\
& \ \mathcal{K}_{\textrm{o}1} \times\cdots \times \mathcal{K}_{\textrm{o}N}\times \mathcal{K}_{\textrm{s}1}\times \cdots \times \mathcal{K}_{\textrm{s}n_\textrm{v}}.
\end{aligned}
\]
\subsection{Persistent excitation}
\label{appendix:PE}
We here provide a detailed discussion of the persistent excitation condition. We consider a linear time-invariant system,
\begin{equation}
\label{eqn:LTI-A}
\begin{aligned}
x(k+1) & = A_1 x(k) + B_1 u(k), \\
y(k) & = C_1 x(k) + D_1 u(k),
\end{aligned}
\end{equation}
where $x \in \mathbb{R}^{\tilde{n}}$, $u \in \mathbb{R}^{\tilde{m}}$, and $y \in \mathbb{R}^{\tilde{p}}$ are the state, input, and output of the system, respectively. Given a length-$T$ pre-collected input-output trajectory $u_\D \in \mathbb{R}^{\tilde{m}T}, y_\D \in \mathbb{R}^{\tilde{p}T}$ and assume the corresponding trajectory of the state is $x_\D \in \mathbb{R}^{\tilde{n}T}$. We can form a length-$L$ trajectory library of the LTI system \eqref{eqn:LTI-A} with the Hankel matrix $H_\D := \col(\mathcal{H}_L(u_\D), \mathcal{H}_L(y_\D))$ whose column vectors are trajectories of the LTI system. The key idea for constructing a data-driven representation of \eqref{eqn:LTI-A} is that $H_\D g$ where $g$ in $\mathbb{R}^{T-L+1}$ can fully capture its behavior with sufficient data. 

It is obvious that $H_\D g$ (\emph{i.e.}, the linear combination of column vectors of $H_\D$) is a valid trajectory of the LTI system due to the linear property. We here mainly illustrate the requirement for the richness of the trajectory library $H_\D$ such that $\mathrm{span}(H_\D)$ can fully represent the trajectory space of the LTI system. It is standard to derive the trajectory library and a valid trajectory of the LTI system as 
\[
H_\D = \begin{bmatrix}
        I_{\tilde{m}L} & \mathbb{0}_{\tilde{m}L \times \tilde{n}} \\
        \mathcal{T}_L & \mathcal{O}_L
    \end{bmatrix} H_0, \quad \begin{bmatrix}
        u \\ y
    \end{bmatrix} = \begin{bmatrix}
        I_{\tilde{m}L} & \mathbb{0}_{\tilde{m}L \times \tilde{n}} \\
        \mathcal{T}_L & \mathcal{O}_L
    \end{bmatrix} \begin{bmatrix}
        x_0 \\  u
    \end{bmatrix}
\]
where we have $H_0 := \col(\mathcal{H}_1(x_{\D, 1:T-L+1}), \mathcal{H}_L(u_{\D}))$ and 
\[
    \mathcal{T}_L = 
    \begin{bmatrix}
    D & 0 & \cdots & 0 \\
    CB & D  & \cdots & 0 \\
    CAB & CB  & \cdots & 0\\
    \vdots  & \vdots & \ddots & \vdots \\
    CA^{L-2}B & CA^{L-3}B  & \cdots & D \\ 
    \end{bmatrix}, \quad  
    \mathcal{O}_L = 
    \begin{bmatrix}
        C\\
        CA\\
        \vdots\\
        CA^{L-1}
    \end{bmatrix}.
\]
For arbitrary $x_0 \in \mathbb{R}^{\tilde{n}}, u \in \mathbb{R}^{\tilde{m}L}$, there exists $g \in \mathbb{R}^{T-L+1}$ such that $\col(x_0, u) = H_0 g$, as long as $H_0$ has full row rank. Thus, the requirement of the richness of the trajectory library is equivalent to making $H_0$ full row rank. It is further proved in \cite{van2020willems} that the persistently exciting of $u_\D$ is a sufficient condition for $H_0$ to be full row rank for a controllable LTI system.

The key challenge in offline data-collection for $\method{dDeeP-LCC}$ is to ensure the input for each CF-LCC subsystem is persistently exciting. These inputs consist of the control signals for the leading CAV of the subsystem as well as its preceding vehicle. 
\begin{itemize}
    \item \textbf{Control input of the CAV:} For the offline data-collection of the CAV, we require a pre-designed controller to maintain its normal driving behavior and some noise can be added to enrich the input signal. In our experiments, we utilize the optimal velocity model (OVM) to collect the input for the CAV with a certain level of i.i.d. noise which is defined as 
    \begin{equation}
    \label{eqn:inputCAV}
    u(t) = \alpha(v_{\des}(s(t))-v(t)) + \beta \dot{s}(t) + \delta_u(t)
     \end{equation}
    where $\delta_u(t) \in \begin{bmatrix}
        -1, 1
    \end{bmatrix}\,\mathrm{m/s^2}$ 
    % and the parameters are set as $\alpha = 0.6$ and $\beta = 0.9$.
    \item \textbf{External input from the preceding vehicle:} The external input for the CF-LCC subsystem comes from the velocity error of its preceding vehicle. In practice, the preceding vehicle is under human control so that its velocity will oscillate around the desired velocity. In our experiment, when the preceding vehicle is the head vehicle, we assume its velocity error is given by
    \begin{equation}
    \label{eqn:inputEx}
    \epsilon(t) = \delta_\epsilon \sim \mathbb{U}\begin{bmatrix}
        -1, 1
    \end{bmatrix}\,\mathrm{m/s}.
    \end{equation}
    For the case that the preceding vehicle is the HDV other than the head vehicle, the velocity error depends on its OVM model which also oscillates around the equilibrium velocity.
\end{itemize}
Thanks to the stochastic nature of the control inputs in \eqref{eqn:inputCAV} and \eqref{eqn:inputEx}, the persistent excitation requirement can be satisfied with a long enough input signal. A longer trajectory results in more columns in the Hankel matrix, increasing the likelihood that it will be full row rank. Moreover, the full rank requirement of the pre-collected data set can be verified before constructing the data-driven representation with it.

\end{document}